\documentclass[preprint,11pt]{elsarticle}
\usepackage{a4}
\usepackage{epsfig}
\usepackage[T1]{fontenc}
\usepackage[utf8]{inputenc}

\usepackage{graphicx,color}

\usepackage{moreverb}

\usepackage{placeins}

\usepackage{amssymb}

\usepackage{amsfonts,amssymb,stmaryrd,amsmath,amssymb}

\newcommand \bn {\boldsymbol{\mathrm{n}}}

\newcommand \bu {\boldsymbol{\mathrm{u}}}
\newcommand \bv {\boldsymbol{\mathrm{v}}}

\newcommand \blambda {\boldsymbol{\mathrm{\lambda}}}
\newcommand \bmu {\boldsymbol{\mathrm{\mu}}}

\newcommand \p {\partial}

\newcommand \x {\mathrm{x}}

\newcommand \R {\mathbb{R}}

\renewcommand \L {\mathrm{L}}

\renewcommand \d {\mathrm{d}}

\renewcommand \div {\mathrm{div}}

\usepackage{fancyhdr}
\fancypagestyle{firststyle}
{
   \fancyhf{}
   \fancyhead[C]{Journal of Fluids and Structures <<Unsteady separation in fluid-structure interaction>>.}
}
               
\journal{Journal of Fluids and Structures}

\begingroup\makeatletter\ifx\SetFigFont\undefined%
\gdef\SetFigFont#1#2#3#4#5{%
\reset@font\fontsize{#1}{#2pt}%
\fontfamily{#3}\fontseries{#4}\fontshape{#5}%
\selectfont}%
\fi\endgroup%

\begin{document}

\begin{frontmatter}


\title{A fictitious domain finite element method for simulations of fluid-structure interactions: The Navier-Stokes equations coupled with a moving solid}
\author{S\'ebastien Court$^*$, Michel Fourni\'e$^{**}$\vspace{0.3cm} \footnote{michel.fournie@math.univ-toulouse.fr}}
\address{$^{*}$Laboratoire de Mathématiques, Campus des Céreaux,\\
Université Blaise Pascal, B.P. 80026, 63171 Aubière cedex, France.\\
$^{**}$Institut de Mathématiques de Toulouse, Unit\'e Mixte C.N.R.S. 5219,\\ 
Universit\'e Paul Sabatier Toulouse III, 118 route de Narbonne, 31062 Toulouse Cedex 9, France.}


\begin{abstract}
The paper extends a stabilized fictitious domain finite element method initially developed for the Stokes problem to the incompressible Navier-Stokes equations coupled with a moving solid. This method presents the advantage to predict an optimal approximation of the normal stress tensor at the interface. The dynamics of the solid is governed by the Newton's laws and the interface between the fluid and the structure is materialized by a level-set which cuts the elements of the mesh. An algorithm is proposed in order to treat the time evolution of the geometry and numerical results are presented on a classical benchmark of the motion of a disk falling in a channel.
\end{abstract}


\end{frontmatter}

\section{Introduction}
Fluid-structure interactions problems remain a challenge both for a comprehensive study of such problems as for the development of robust numerical methods (see a review in \cite{Hou&Wang&Layton}). One class of numerical methods is based on meshes that are conformed to the interface where the physical boundary conditions are imposed \cite{LT, SMSTT0, SST}. As the geometry of the fluid domain changes through the time, re-meshing is needed, which is excessively time-consuming, in particular for complex systems. An other class of numerical methods is based on non-conforming mesh with a fictitious domain approach where the mesh is cut by the boundary. Most of the non-conforming mesh methods are based on the immersed boundary methods where force-equivalent terms are added to the fluid equations in order to represent the fluid structure interaction \cite{Peskinacta, Mittal&Iaccarino}. Many related numerical methods have been developed, in particular the popular distributed Lagrange multiplier method, introduced for rigid bodies moving in an incompressible flow \cite{Glowinski}. In this method, the fluid domain is extended in order to cover the rigid domain where the fluid velocity is required to be equal to the rigid body velocity.\\
More recently, eXtended Finite Element Method introduced by Mo\"{e}s, Dolbow and Belytschko in \cite{MoesD} (see a review of such methods in \cite{reviewXfem}) has been adapted to fluid structure interactions problems in \cite{MoesB, SukumarC, Gerstenberger2008, Choi2010}. The idea is similar to the fictitious domain / Lagrange multiplier method aforementioned, but the fluid velocity is no longer extended inside the structure domain, and its value given by the structure velocity is enforced by a Lagrange multiplier only on the fluid-structure interface. One thus gets rid of unnecessary fluid unknowns. Besides, one easily recovers the normal trace of the Cauchy stress tensor on the interface. We note that this  method has been originally developed for problems in structural mechanics mostly in the context of cracked domains, see for example \cite{HaslR, MoesG, Stazi, SukumarM, Stolarska}. The specificity of the method is that it combines a level-set representation of the geometry of the crack with an enrichment of a finite element space by singular and discontinuous functions.\\
In the context of fluid-structure interactions, the difficulty related to the applications of such techniques lies in the choice of the Lagrange multiplier space used in order to take into account the interface, which is not trivial because of the fact that the interface cuts the mesh (see \cite{Bechet2009} for instance). 
In particular, the natural mesh given by the points of intersection of the interface with the global mesh cannot be used directly. An algorithm to construct a multiplier space satisfying the inf-sup condition is developed in \cite{Bechet2009}, but its implementation can be difficult in practice. 
The method proposed in the present paper tackles this difficulty by using a stabilization technique proposed in \cite{HaslR}. This method was adapted to contact problems in elastostatics in \cite{HildR} and more recently to the Stokes problem in \cite{CourtFournieLozinski}. 
An important feature of this method (based on the eXtended Finite Element Method approach, similarly to \cite{Gerstenberger2008, Choi2010}) is that the Lagrange multiplier is identified with the normal trace of the Cauchy stress tensor $\sigma(\bu,p)\bn$ at the interface.  
Moreover, it is possible to obtain a good numerical approximation of $\sigma(\bu,p)\bn$  (the proof is given in \cite{CourtFournieLozinski} for the Stokes problem). This property is crucial in fluid-structure interactions since this quantity gives the force exerted by the viscous fluid on the structure. In the present paper, we propose to extend this method to the Navier-Stokes equations coupled with a moving solid. Note that alternative methods based on the Nitsche's work \cite{Nitsche} (such as \cite{Burman1, Burman3} in the context of the Poisson problem and \cite{Massing} in the context of the Stokes problems) do not introduce the Lagrange multiplier and thus do not necessarily provide a good numerical approximation of this force. Our method based on boundary forces is particular interesting for control flow around a structure. The control function can be localized on the boundary of the structure where we impose its local deformation. In order to perform direct numerical simulations of such a control, efficient tools based on accurate computations on the interface must be developed. The present approach is one brick in this research topic where recent development towards stabilized Navier-Stokes equations are proposed (like in \cite{JPR}).\\

The outline of the paper is as follows. The continuous fluid-structure interactions problem is given in Section~\ref{section2} and the weak formulation with the introduction of a Lagrange multiplier for imposing the boundary condition at the interface is given in Section~\ref{subsection2}. Next, in Section~\ref{section3} the fictitious domain method is recalled with the introduction of the finite element method (Section~\ref{fem}) with a time discretization (Section~\ref{time}). Section~\ref{section4} is devoted to numerical tests and validation on a benchmark corresponding to the falling of a disk in a channel. The efficiency of the method is presented before conclusion.

\section{The model} \label{section2}
\subsection{Fluid-structure interactions}
\hspace*{0.5cm} We consider a moving solid which occupies {a time-depending domain denoted by $\mathcal{S}(t)$}. The remaining domain $\mathcal{F}(t) = \mathcal{O} \setminus \overline{\mathcal{S}(t)}$ corresponds to the fluid flow.

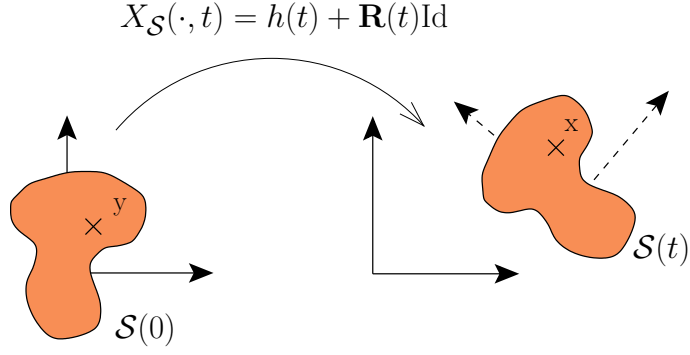
\begin{figure}[!h]
\begin{center}
\scalebox{0.5}{
\includegraphics[trim = 0cm 0cm 0cm 0cm, clip]{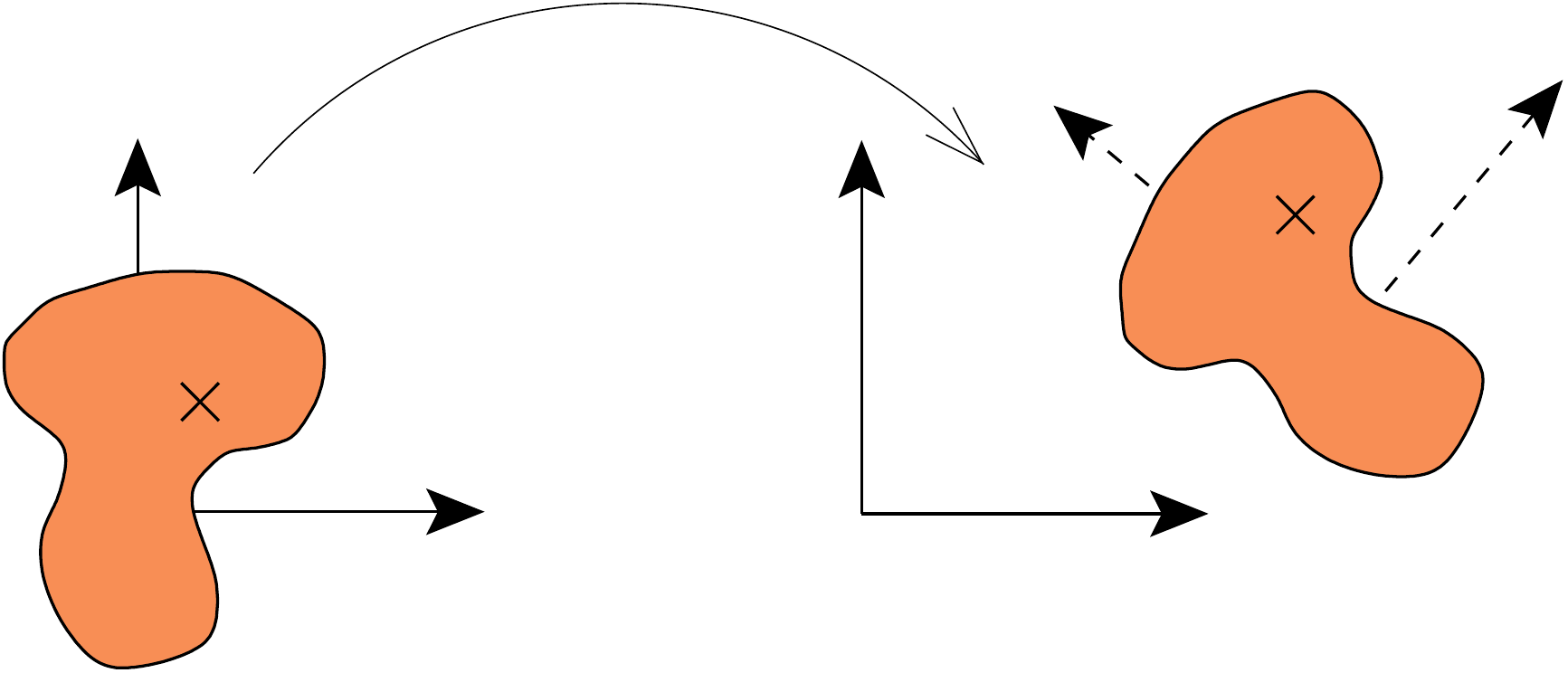}
\setlength{\unitlength}{4144sp}
\begin{picture}(0,0)
\put(-6706,3750){\makebox(0,0)[lb]{\smash{{\color[rgb]{0,0,0}\begin{Huge}$X_{\mathcal{S}}(\cdot,t)
= h(t) + \mathbf{R}(t)\mathrm{Id}$\end{Huge}}%
}}}
\put(-1386,2516){\makebox(0,0)[lb]{\smash{{\color[rgb]{0,0,0}\begin{huge}$\mathrm{x}$\end{huge}}%
}}}
\put(-6786,1561){\makebox(0,0)[lb]{\smash{{\color[rgb]{0,0,0}\begin{huge}$\mathrm{y}$\end{huge}}%
}}}
\put(-6766,30){\makebox(0,0)[lb]{\smash{{\color[rgb]{0,0,0}\begin{Huge}$\mathcal{S}(0)$\end{Huge}}%
}}}
\put(-566,986){\makebox(0,0)[lb]{\smash{{\color[rgb]{0,0,0}\begin{Huge}$\mathcal{S}(t)$\end{Huge}}%
}}}
\end{picture}%
}
\end{center}
\centering \caption{Decomposition of the solid movement.\label{fig1}}
\end{figure}
\FloatBarrier

The displacement of a rigid solid can be given by the knowledge of ${\bf h}(t)$, namely the position of its gravity center, and $\mathbf{R}(t)$ its rotation given by $\displaystyle \left( \begin{array}{cc} c& -s \\ s& c\\ \end{array} \right)$ for $\displaystyle c=\cos(\theta(t))$, $s=\sin(\theta(t))$, where $\theta(t)$ is the rotation angle of the solid (see Figure~\ref{fig1}).
Then at time $t$ the domain occupied by the structure is given by 
\begin{eqnarray*}
\displaystyle \mathcal{S}(t) & = & {\bf h}(t) + \mathbf{R}(t)\mathcal{S}(0).
\end{eqnarray*}
{\it Remark:} This formulation can be extended in order to consider general deformations of the structure. Then we would have to define a mapping $X^*(\cdot,t)$ which corresponds to the deformation of the solid in its own frame of reference. Then, $\mathcal{S}(t) = X_S(\mathcal{S}(0),t)$ where $X_S(\mathrm{y},t)  = {\bf h}(t) + \mathbf{R}(t)X^*(\mathrm{y},t)$, for $\mathrm{y} \in \mathcal{S}(0)$.\\

The velocity of the incompressible viscous fluid of density $\rho_f$ is denoted by ${\bf u}$, the pressure by $p$ and $\nu$ is the dynamic viscosity. We denote by $\bn$ the outward unit normal vector to $\p \mathcal{F}$ (the boundary of $\mathcal{F}$), and the normal trace on the interface $\Gamma = \partial \mathcal{S}(t)$ of the Cauchy stress tensor is given by
\begin{eqnarray*}
\sigma(\bu,p)\bn  =   2\nu D(\bu){\bf n} -p\bn, & & \text{ with } D(\bu)  =  \frac{1}{2} \left(\nabla \bu + \nabla \bu^T \right).
\end{eqnarray*}
When gravity forces are considered (we denote by ${\bf g}$ the gravity field), the fluid flow is modeled by the incompressible Navier-Stokes equations

\begin{equation}
\label{eq1}
\left \{
\begin{array}{llr}
\displaystyle \rho_f\left ( \frac{\partial {\bf u}}{\partial t} + ({\bf u}.\nabla){\bf u}  \right ) - \nu \Delta {\bf u} +\nabla p = \rho_f {\bf g}, & \mathrm{x} \in \mathcal{F}(t), & t \in (0,T),\\
\displaystyle \mbox{div}({\bf u}) =0, & \mathrm{x} \in \mathcal{F}(t), & t \in (0,T),\\
\displaystyle {\bf u} = 0, & \mathrm{x} \in \partial \mathcal{O}  , &  t \in (0,T),\\
\end{array}
\right .
\end{equation}
and the Newton's laws are considered for the dynamics of the solid
\begin{equation}
\label{eq2}
\left \{
\begin{array}{l}
\displaystyle m_s{\bf h}''(t) = - \int_{\partial \mathcal{S}(t)} {\sigma({\bf u},p){\bf n}} d \Gamma -m_s {\bf g},\\
\displaystyle I\theta''(t) =  - \int_{\partial \mathcal{S}(t)}({\bf x}-{\bf h}(t))^{\bot}\cdot {\sigma({\bu},p){\bf n}} d \Gamma,
\end{array}
\right .
\end{equation}
where $m_s$ is the mass of the solid, and $I$ is its moment of inertia.\\
At the interface $\partial \mathcal{S}(t)$, for the coupling between fluid and structure, we impose the continuity of the velocity 
\begin{equation}
\label{eq3}
\displaystyle {\bf u}({\bf x},t)  =  {\bf h}'(t) + \theta'(t) ({\bf x}-{\bf h}(t))^{\bot}={\bf u_{\Gamma}}, \ \ \  {\bf x} \in \partial \mathcal{S}(t), \ \ \   t \in (0,T).
\end{equation}
The coupled system (\ref{eq1})--(\ref{eq3}) has for unknowns ${\bf u}$, $p$, ${\bf h}(t)$ and the angular velocity $\omega(t) = \theta'(t) $ (a scalar function in 2D).

\subsection{Weak formulation of the problem with stabilization terms} \label{subsection2}
We consider the coupled system (\ref{eq1})-(\ref{eq3}) and we assume that the boundary condition imposed at the interface $\Gamma = \partial \mathcal{S}(t)$ is sufficiently regular to make sense, and we introduce the following functional spaces (based on the classical Sobolev spaces $\L^2(\mathcal{F})$, $\mathbf{H}^1(\mathcal{F})$, $\mathbf{H}^{1/2}(\Gamma)$ and $\mathbf{H}^{-1/2}(\Gamma)$, see \cite{Evans} for instance)
\begin{eqnarray*}
\begin{array}{lcl}
\mathbf{V} &= &\left\{ \bv\in \mathbf{H}^1(\mathcal{F}) \mid \bv=0 \text{ on } \p \mathcal{O} \right\},\\
Q &=& \L^2_0(\mathcal{F}) = \left\{p \in \L^2(\mathcal{F}) \mid \displaystyle \int_{ \mathcal{F}}p\ \d \mathcal{F} = 0 \right\},\\
\mathbf{W} &=& \mathbf{H}^{-1/2}(\Gamma) = \left(\mathbf{H}^{1/2}(\Gamma) \right)'.\\
\end{array}
\end{eqnarray*}
Due to the fact that we only consider boundary conditions of Dirichlet type, we impose to the pressure $p$ to have null average (this condition is taken into account in $Q$). That variational formulation can be done in three steps:
\begin{itemize}
\item[Step 1 --] Classical formulation of the Navier-Stokes and structure equations (in the formulation $\gamma$ and $\blambda$ are equal to $0$);
\item[Step 2 --] Introduction of Lagrange multiplier $\blambda$ in order to take into account the Dirichlet condition at the interface $\Gamma$ (in the formulation only $\gamma$ is equal to $0$);
\item[Step 3 --] Introduction of stabilization terms with a parameter $\gamma$.
\end{itemize}

This new unknown $\blambda$ plays a critical role due to the fact that it is equal to the 
normal trace of the Cauchy stress tensor $\sigma(\bu,p)\bn$ (this equality is described in \cite{Gunzburger}). The stabilization terms are associated with the constant parameter $\gamma$ (chosen sufficiently small). The variational problem that we consider is the following:
\begin{eqnarray*}
\hspace*{-15pt} & & \text{Find $({\bf u},p,{\bf \lambda},{\bf h}',{\bf h},\theta',\theta) \in \mathbf{V} \times Q \times \mathbf{W} \times \R^2   \times \R^2  \times \R \times \R$ such that} \\
\hspace*{-15pt} & & \left\{ \begin{array} {ll}
\displaystyle \int_{\mathcal{F} } \rho_f \frac{\partial {\bf u}}{\partial t}\cdot{\bf v} \mathrm{d}\mathcal{F} +  \mathcal{A}(({\bf u},p,{\blambda});{\bf v}) + \int_{\mathcal{F}}  \rho_f [({\bf u}\cdot\nabla){\bf u}]\cdot {\bf v} \mathrm{d}\mathcal{F}  =  \int_{\mathcal{F} } \rho_f  {\bf g} \cdot {\bf v}\mathrm{d}\mathcal{F},
 & \forall {\bf v} \in \mathbf{V}, \\
\mathcal{B}(({\bf u},p,\boldsymbol{\lambda});q)  = 0,  & \forall q \in Q, \\
\mathcal{C}(({\bf u},p,\boldsymbol{\lambda});{\bmu}) = \mathcal{G}({\bmu}),  \quad & \forall {\bf \bmu} \in \mathbf{W},\\
\displaystyle m_s{\bf h}''(t) = - \int_{\partial \mathcal{S}(t)}  {\blambda} d \Gamma -m_s {\bf g},&\\

\displaystyle I\theta''(t) =  - \int_{\partial \mathcal{S}(t)}({\bf x}-{\bf h}(t))^{\bot}\cdot \boldsymbol{\lambda} d \Gamma,&
\end{array} \right. \label{FVaugmented}
\end{eqnarray*}
where
\begin{eqnarray*}
\mathcal{A}(({\bf u},p,\boldsymbol{\lambda});{\bf v}) & = & 2\nu\int_{\mathcal{F}}D({\bf u}):D({\bf v})\mathrm{d}\mathcal{F} - \int_{\mathcal{F}}p\div \ {\bf v}\mathrm{d}\mathcal{F}  - \int_{\Gamma} \boldsymbol{\lambda}\cdot {\bf v}\mathrm{d}\Gamma \\
& &\hspace*{-2.5cm}    -4\nu^2\gamma \int_{\Gamma}\left( D({\bf u})\bn \right)\cdot \left( D({\bf v})\bn \right)\mathrm{d}\Gamma +2\nu \gamma \int_{\Gamma}p \left(D({\bf v})\bn\cdot \bn \right)\mathrm{d}\Gamma +2\nu \gamma \int_{\Gamma} \boldsymbol{\lambda} \cdot \left(D({\bf v})\bn\right)\mathrm{d}\Gamma , \\
\mathcal{B}(({\bf u},p,\boldsymbol{\lambda});q) & = & - \int_{\mathcal{F}}q\div\ {\bf u}\mathrm{d}\mathcal{F} +2\nu \gamma \int_{\Gamma}q\left(D({\bf u})\bn\cdot \bn \right)\mathrm{d}\Gamma -\gamma \int_{\Gamma}pq \mathrm{d}\Gamma - \gamma \int_{\Gamma} q\boldsymbol{\lambda} \cdot \bn\mathrm{d}\Gamma , \\
\mathcal{C}(({\bf u},p,\boldsymbol{\lambda});\boldsymbol{\mu}) & = & -\int_{\Gamma} \boldsymbol{\mu} \cdot {\bf u}\mathrm{d}\Gamma +2\nu \gamma \int_{\Gamma}\boldsymbol{\mu} \cdot (D({\bf u})\bn)\mathrm{d}\Gamma -\gamma \int_{\Gamma}p(\boldsymbol{\mu}\cdot \bn)\mathrm{d}\Gamma - \gamma \int_{\Gamma} \boldsymbol{\lambda} \cdot \boldsymbol{\mu} \mathrm{d}\Gamma,\\
\hspace*{-0.8cm} \mathcal{G}(\boldsymbol{\mu})  &= & -\int_{\Gamma} \boldsymbol{\mu} \cdot {\bf u_{\Gamma}} \mathrm{d}\Gamma = - \int_{\Gamma} \boldsymbol{\mu} \cdot  ({\bf h}'(t) + \theta'(t) ({\bf x}-{\bf h}(t))^{\bot})\mathrm{d}\Gamma.\\
\end{eqnarray*}
{\it Remark:} The formulation can be justified by the introduction of an extended Lagrangian - \`a la Barbosa-Hughes, see \cite{Barbosa1} - whose a stationary point is a weak solution of the problem. The first-order derivatives of this Lagrangian leads to forcing $\blambda$ to reach the desired value corresponding to $\sigma(\bu,p)\bn$.

\section{Fictitious domain approach}\label{section3}
We refer to the article \cite{CourtFournieLozinski} for the details of the fictitious domain approach we consider here. In the following, we recall the method used for the present work.

\subsection{Finite element discretization} \label{fem}
The fictitious domain for the fluid is considered on the whole domain  $\mathcal{O}$. Let us introduce three discrete finite element spaces, {$\tilde{\mathbf{V}}^h \subset \mathbf{H}^1(\mathcal{O})$, $\tilde{Q}^h \subset \L^2_0(\mathcal{O})$ and $\tilde{\mathbf{W}}^h \subset \mathbf{L}^2(\mathcal{O})$}. Notice that the spaces $\mathbf{V},\  Q,\  \mathbf{W}$  introduced to define the weak formulation are included into those spaces defined all over the domain $\mathcal{O} = \mathcal{F} \cup  \mathcal{S} $. In practice, $\mathcal{O}$ is a simple domain, so that the construction of a unique mesh for all spaces is straightforward (the interface between the fluid and the structure is not considered). Let us consider for instance a rectangular domain where a structured uniform mesh $\mathcal{T}^h$ can be constructed (see Figure~\ref{figMesh}). Classical finite element discretizations can be defined on the spaces {$\tilde{\mathbf{V}}^h$, $\tilde{Q}^h$ and $\tilde{\mathbf{W}}^h$}. For $\tilde{\mathbf{V}}^h$, let us consider for instance a subspace of the continuous functions $ C(\overline{\mathcal{O}})$ defined by
\begin{eqnarray*}
\tilde{\mathbf{V}}^h & = & \left\{\bv^h \in C(\overline{\mathcal{O}})\mid \bv^h_{\left| \p \mathcal{O}\right.} = 0, \  \bv^h_{\left| T\right.} \in P(T), \ \forall T \in \mathcal{T}^h \right\}, \label{defvtilde}
\end{eqnarray*}
where $P(T)$ is a finite dimensional space of regular functions, containing $P_k(T)$ the polynom space of degree less or equal to an integer $k$ ($k \geq 1$). For more details, see \cite{Ern} for instance. The mesh step stands for $\displaystyle h = \max_{T\in \mathcal{T}^h} h_T$, where $h_T$ is the diameter of $T$. In order to split the fluid domain and the structure domain, we define spaces on the fluid part $\mathcal{F}$ and on the interface $\Gamma$ only, as 
\begin{eqnarray*}
\mathbf{V}^h := \tilde{\mathbf{V}}^h_{\left| \mathcal{F} \right.}, \quad Q^h := \tilde{Q}^h_{\left|\mathcal{F}\right.}, \quad  \mathbf{W}^h := \tilde{\mathbf{W}}^h_{\left| \Gamma \right.}.
\end{eqnarray*}

\begin{center}
\begin{figure}[h!]
\includegraphics[trim = 0cm 0cm 0cm 0cm, clip, scale= 0.50]{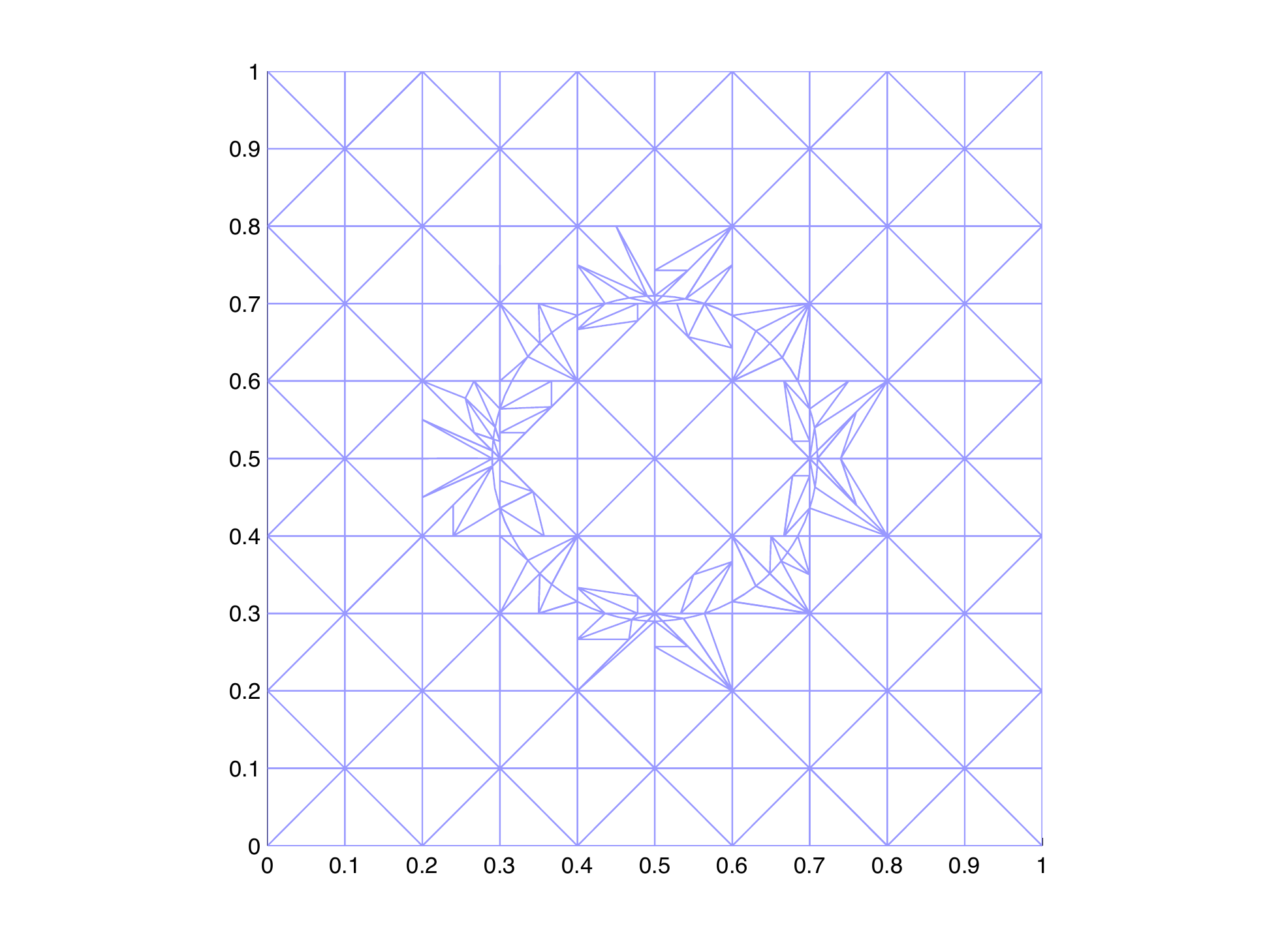}
\centering \caption{Illustration of the elements cut with respect to the level-set.\label{figsupercut}}
\end{figure}
\end{center}
\FloatBarrier

Notice that $\mathbf{V}^h $, $Q^h$, $\mathbf{W}^h$ are respective natural discretizations of $\mathbf{V}$, $Q$ and $\mathbf{W}$. It corresponds to cutting the basis functions of spaces $\tilde{\mathbf{V}}^h$, $Q^h$ and $\tilde{\mathbf{W}}^h$, as shown in Figure~\ref{figsupercut}. This approach is equivalent to the eXtended Finite Element Method, as proposed in \cite{Choi2010} or \cite{Gerstenberger2008}, where the standard finite element method basis functions are multiplied by Heaviside functions ($H({\bf x}) = 1$ for ${\bf x} \in \mathcal{F}$ and $H({\bf x})=0$ for ${\bf x}\in \mathcal{O} \setminus \mathcal{F}$), and the products are substituted in the variational formulation of the problem. Thus the degrees of freedom inside the fluid domain $\mathcal{F}$ are used in the same way as in the standard finite element method, whereas the degrees of freedom in the solid domain $\mathcal{S}$ at the vertexes of the elements cut by the interface (the so called virtual degrees of freedom) do not define the field variable at these nodes, but they are necessary to define the fields on $\mathcal{F}$ and to compute the integrals over $\mathcal{F}$. The remaining degrees of freedom, corresponding to the basis functions with support completely outside of the fluid, are eliminated (see Figure~\ref{figMesh}). We refer to the papers aforementioned for more details.

\begin{figure}[h!]
\begin{center}
\hspace*{1cm} \includegraphics[trim = 9cm 2cm 2cm 1cm, clip, scale= 0.50]{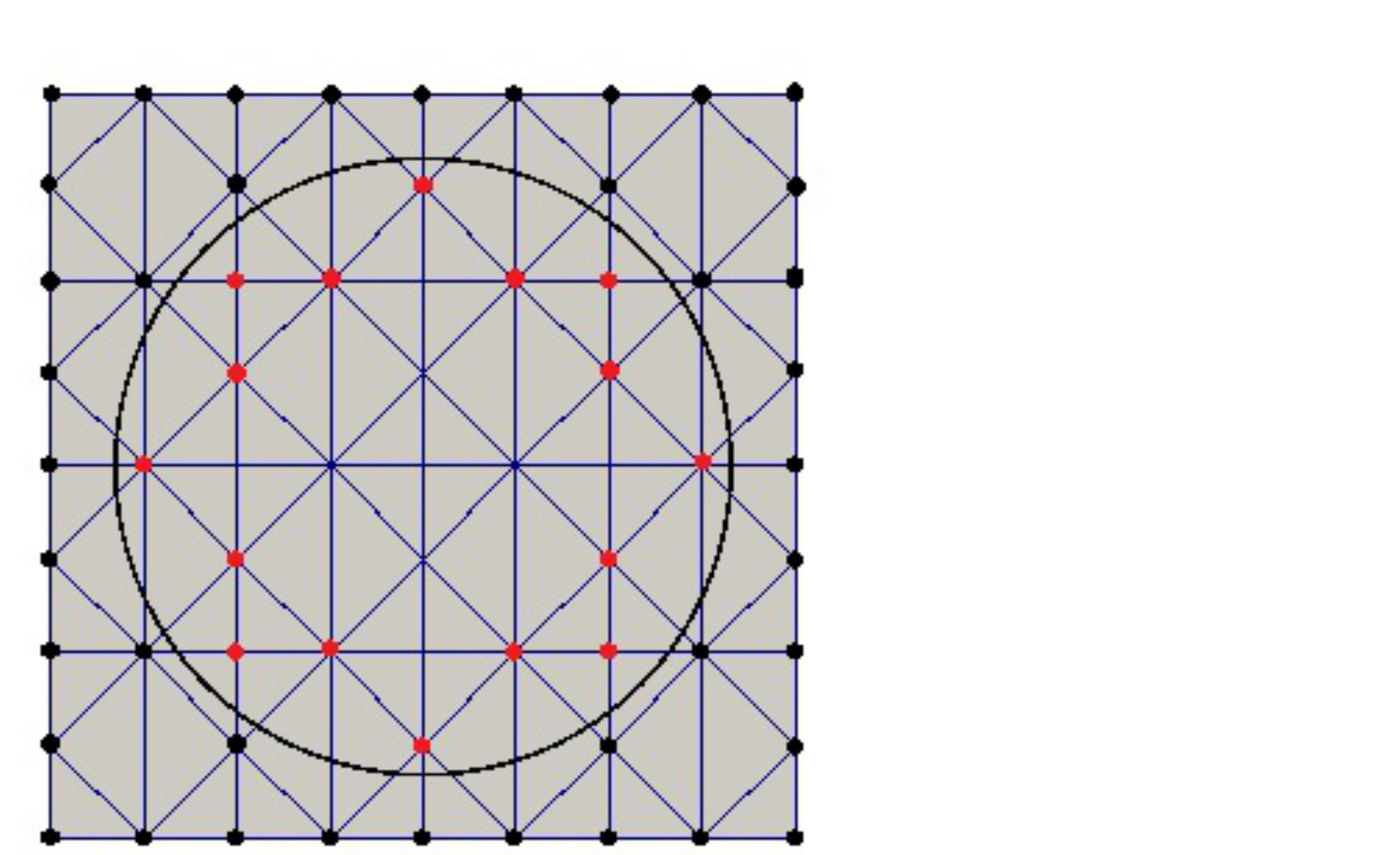}\hspace*{-2cm}  
\includegraphics[trim = 9cm 2cm 2cm 1cm, clip, scale= 0.50]{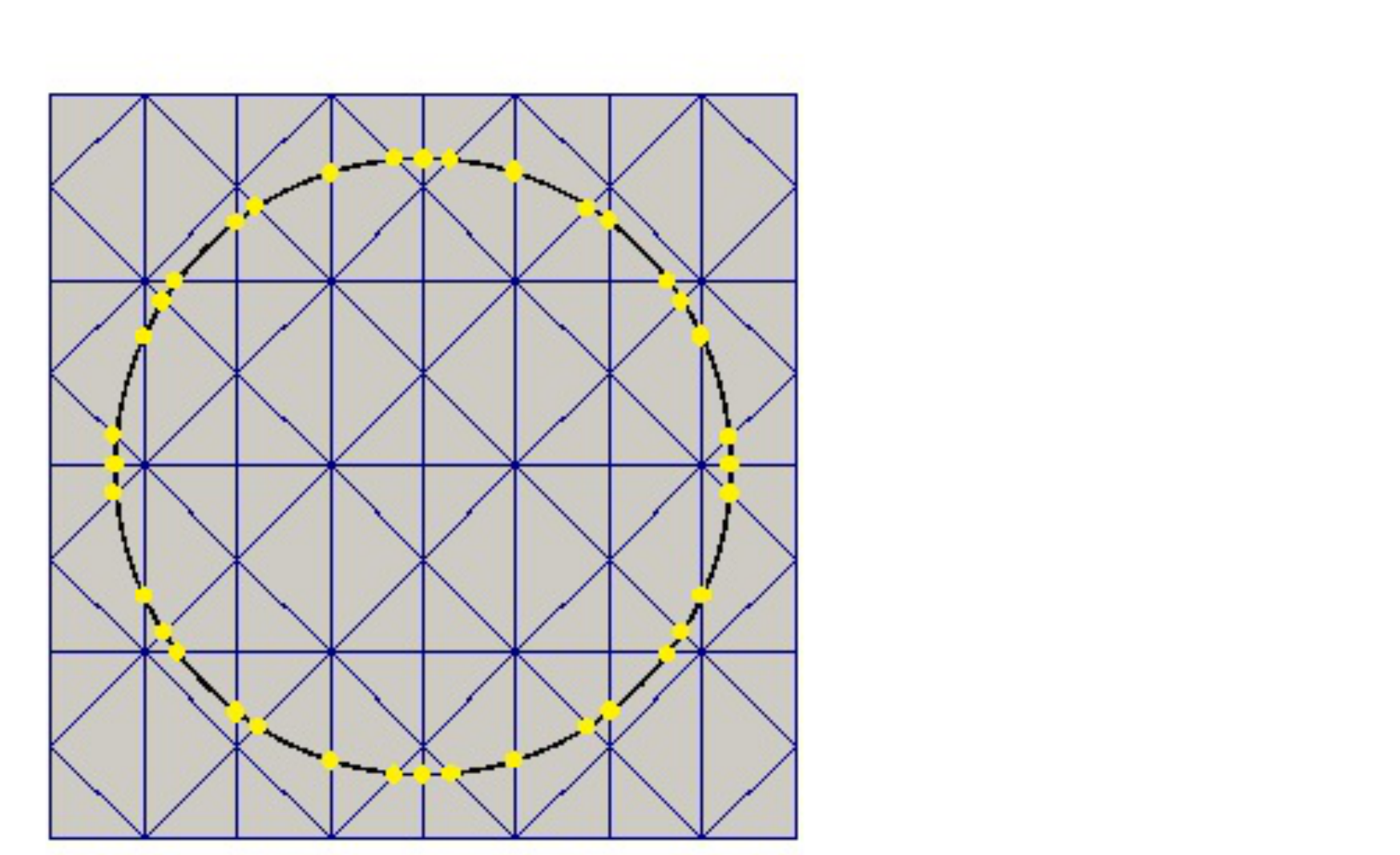}

\hspace*{-0.5cm} (a) \hspace*{6.5cm} (b)\\
\caption{\label{figMesh} Mesh on a fictitious domain.
(a) Standard degrees of freedom (black, outside the disk), virtual ones (red, inside the disk), remaining ones are removed. (b) Bases nodes used for the multiplier space (yellow, on the disk).}
\end{center}
\end{figure}
\FloatBarrier

The discrete problem consists in finding $(\bu^h,p^h,\blambda^h,{\bf h}',{\bf h},\theta', \theta) \in \mathbf{V}^h \times Q^h \times \mathbf{W}^h \times \R^2  \times \R^2 \times \R \times \R$ such that 
\begin{eqnarray*}
 \left\{ \begin{array} {lll}
\displaystyle \int_{\mathcal{F} } \rho_f \frac{\partial {\bf u}^h}{\partial t}\cdot{\bf v^h} \mathrm{d}\mathcal{F} +  \mathcal{A}(({\bf u}^h,p^h,{\blambda}^h);{\bf v}^h)  +\int_{\mathcal{F}} \rho_f [({\bf u}^h\cdot\nabla){\bf u}^h]\cdot{\bf v}^h \mathrm{d}\mathcal{F}  =  \int_{\mathcal{F} } \rho_f  {\bf g} \cdot{\bf v}^h\mathrm{d}\mathcal{F},\\
&   \hspace*{-2cm} \forall {\bf v}^h \in \mathbf{V}^h, \\ &\\
\mathcal{B}(({\bf u}^h,p^h,{\blambda}^h);q^h)  = 0, & \hspace*{-2cm}  \forall q^h \in Q^h, \\
&\\
\mathcal{C}(({\bf u}^h,p^h,{\blambda}^h);{\bmu}^h) = \mathcal{G}({\bmu}^h), &  \hspace*{-2cm}  \forall {\bmu}^h \in \mathbf{W}^h,\\
&\\
\displaystyle m_s{\bf h}''(t) = - \int_{\partial \mathcal{S}(t)}  {\blambda}^h \d \Gamma -m_s {\bf g},\qquad 
\displaystyle I\theta''(t) =  - \int_{\partial \mathcal{S}(t)}({\bf x}-{\bf h}(t))^{\bot}\cdot {\blambda}^h \d \Gamma.&\\
\end{array} \right. \label{FVhaugmented}
\end{eqnarray*}
This is a system of nonlinear differential algebraic equations which can be formulated into a compact form. We denote by $\boldsymbol{U}$, $\boldsymbol{P}$ and $\boldsymbol{\Lambda}$ the respective degrees of freedom of $\bu^h$, $p^h$ and $\blambda^h$. After standard finite element discretization of the following bilinear forms
\begin{eqnarray*}
\mathcal{M}_{{\bf u}{\bf u}} : ({\bf u},{\bf v}) & \longmapsto & \int_{\mathcal{F}}\rho_f{\bf u}.{\bf v}\mathrm{d}\mathcal{F}, \qquad
\mathcal{M}_{\blambda} : {\blambda}   \longmapsto  -\int_{\partial \mathcal{S}(t)}{\blambda} \mathrm{d}\Gamma,\\
\mathcal{A}_{{\bf u}{\bf u}} : ({\bf u},{\bf v}) & \longmapsto & 2\nu\int_{\mathcal{F}}D({\bf u}):D({\bf v})\mathrm{d}\mathcal{F} - 4\nu^2\gamma \int_{\Gamma}\left( D({\bf u})\bn \right)\cdot \left( D({\bf v})\bn\right)\mathrm{d}\Gamma , \\
\mathcal{A}_{{\bf u}p} : ({\bf v},p) & \longmapsto & - \int_{\mathcal{F}}p\div \ {\bf v} \mathrm{d}\mathcal{F} + 2\nu \gamma \int_{\Gamma}p \left(D({\bf v})\bn\cdot \bn \right)\mathrm{d}\Gamma , \\
\mathcal{A}_{{\bf u}{\blambda}} : ({\bf u},{\blambda}) & \longmapsto & - \int_{\Gamma} {\blambda}\cdot {\bf v}\mathrm{d}\Gamma + 2\nu \gamma \int_{\Gamma} {\blambda} \cdot \left(D({\bf v})\bn\right)\mathrm{d}\Gamma , \\
\mathcal{A}_{pp} : (p,q) & \longmapsto & -\gamma \int_{\Gamma}pq \mathrm{d}\Gamma , \hspace*{0.5cm} 
\mathcal{A}_{p{\blambda}} : (q,{\blambda})  \longmapsto  - \gamma \int_{\Gamma} q{\blambda} \cdot \bn\mathrm{d}\Gamma , \\
\mathcal{A}_{{\blambda}{\blambda}} : ({\blambda},{\bmu}) & \longmapsto & -\gamma \int_{\Gamma} {\blambda} \cdot {\bmu} \mathrm{d}\Gamma,
\end{eqnarray*}
we define matrices like $M_{{\bf u}{\bf u}}$ from  $\mathcal{M}_{{\bf u}{\bf u}}$, etc..., the vector $\boldsymbol{G}$ from $\mathcal{G}$, ${\bf F}$ from the gravity forces $\rho_f {\bf g}$, $M_{\blambda}$ the matrix computed by integration over $\Gamma$ of the $\mathbf{W}^h$ basis functions and $N(\boldsymbol{U}(t))\boldsymbol{U}(t)$ the matrix depending on the velocity and corresponding to the nonlinear convective term $\displaystyle \int_{\mathcal{F}} \rho_f[({\bf u} \cdot \nabla){\bf u}] \cdot{\bf v} \mathrm{d}\mathcal{F}$. Then the matrix formulation is given by
\begin{eqnarray}
M_{{\bf u}{\bf u}} \frac{\mbox{d} \boldsymbol{U(t)}}{\mbox{d} t} +A_{{\bf u}{\bf u}} \boldsymbol{U}(t) + N( \boldsymbol{U}(t)) \boldsymbol{U}(t) +  A_{{\bf u}p}\boldsymbol{P}(t) +  A_{{\bf u}{\bf \lambda}}\boldsymbol{\Lambda}(t)  = {\bf F}, \label{M1} \\
A^T_{{\bf u}p}  \boldsymbol{U}(t) +  A_{pp} \boldsymbol{P}(t)  + A_{p{\bf \lambda}}\boldsymbol{\Lambda}(t) = 0, \label{M2} \\
A^T_{{\bf u}{\bf \lambda}}  \boldsymbol{U}(t)  +  A^T_{p{\bf \lambda}} \boldsymbol{P}(t)  +  A_{{\bf \lambda} {\bf \lambda}}\boldsymbol{\Lambda}(t) = \boldsymbol{G},\label{M3}\\
\displaystyle m_s{\bf h}''(t) = M_{ \bf \lambda} \boldsymbol{\Lambda}(t) -m_s {\bf g},&&\label{M4}\\
\displaystyle \text{{$I\theta''(t) =  M_{ \bf \lambda}
\left[ ({\bf x}-{\bf h}(t))^{\bot} \cdot \boldsymbol{\Lambda}(t)\right]$}}.	&&\label{M5}
\end{eqnarray}
At the interface $\Gamma = \partial \mathcal{S}(t)$ represented by a level-set function which cuts the global mesh, the coupling between the fluid and the structure is imposed by a Dirichlet condition whose elements are determined through the computation of $\sigma(\bu,p){\bf n}$. The main advantage of our numerical method - mathematically justified in \cite{CourtFournieLozinski} - is to return an optimal approximation $\boldsymbol{\Lambda}(t)$ of $\sigma(\bu,p){\bf n}$ at the interface. Getting a good approximation of this quantity is crucial for the dynamics of the system.

\subsection{Time discretization and treatment of the nonlinearity} \label{time}
Classical methods like $\theta$-methods can be used for the time discretization. For a matter of unconditional stability of the scheme, we consider an implicit discretization based on the backward Euler method. We denote by $\boldsymbol{U}^{n+1}$ the solution at the time level $t^{n+1}$ and $dt = t^{n+1}-t^n$ is the time step. Particular attention must be done for a moving particle problem. Indeed, at the time level $t^{n+1}$ the solid occupies $\mathcal{S}(t^{n+1})$ which is different from the previous time level $t^n$. So, the field variable at the time level $t^{n+1}$ can become undefined near the interface since there was no fluid flow at the time level $t^n$ ($\mathcal{S}(t^{n+1})\neq \mathcal{S}(t^{n})$ for the solid and $\mathcal{F}(t^{n+1})\neq \mathcal{F}(t^{n})$ for the fluid). In other words, some degrees of freedom for the fluid part which are not considered at the time level $t^n$ must be taken into account at the time level $t^{n+1}$. In particular, the velocity field must be known in such nodes. In the present work, we impose the velocity to be equal to the motion of the solid. The validity of this approximation is justified as soon as time step is sufficiently small to ensure that the structure moves progressively across the mesh without jump of cells (when the level-set doesn't cuts this cell). This constraint  is not too strong and corresponds to the classical CFL condition for velocity of the structure.\\
In the following we present the algorithm we perform to compute at the time level $t^{n+1}$ the solution ($\boldsymbol{U}^{n+1}, \boldsymbol{P}^{n+1}, \boldsymbol{\Lambda}^{n+1},{\bf h}'^{n+1},{\bf h}^{n+1},\theta'^{n+1},\theta^{n+1}$) on $\mathcal{F}(t^{n+1})$. To simplify, we assume that $dt$ is constant. At the time level $t^n$ we have access to ($\boldsymbol{U}^n, \boldsymbol{P}^n, \boldsymbol{\Lambda}^n,{\bf h}'^n,{\bf h}^n,\theta'^n,\theta^{n}$) on $\mathcal{F}(t^{n})$.
\begin{itemize}
\item[1--] {\bf Velocity of the structure} -From $\boldsymbol{\Lambda}^n$, we compute $({\bf h}'^{n+1},\theta'^{n+1})$ using (\ref{M4}) and (\ref{M5}) with the mid-point method, as
\begin{eqnarray*}
\displaystyle m_s \frac{{\bf h}'^{n+1}-{\bf h}'^n}{dt}  =   M_{ \bf \lambda} \boldsymbol{\Lambda}^{n}  -m_s {\bf g},&&\\
\displaystyle I\frac{\theta'^{n+1} - \theta'^n}{dt} =
M_{ \bf \lambda}  \left [ ({\bf x}-{\bf h}^{n})^{\bot}\cdot  \boldsymbol{\Lambda}^{n}\right ] .&&\\
\end{eqnarray*}
\item[2--] {\bf Position of the structure} - From $\boldsymbol{\Lambda}^n$, we compute $({\bf h}^{n+1},\theta^{n+1})$ using (\ref{M4}) and (\ref{M5}) with  mid point rule
\begin{eqnarray*}
\displaystyle m_s \frac{{\bf h}^{n+1}-2{\bf h}^n+{\bf h}^{n-1}}{dt^2}  =  M_{ \bf \lambda} \boldsymbol{\Lambda}^{n} -m_s {\bf g},&&\\
\displaystyle I\frac{\theta^{n+1} - 2 \theta^n +  \theta^{n-1}}{dt^2} = M_{ \bf \lambda}  \left [ ({\bf x}-{\bf h}^{n})^{\bot}\cdot  \boldsymbol{\Lambda}^{n}\right ] .&&
\end{eqnarray*}
\item[3--] We update the geometry to determine $\mathcal{F}(t^{n+1})$. It corresponds to update the position of the level-set which is defined from ${\bf h}^{n+1}$ and $\theta^{n+1}$.
\item[4--] We complete the velocity $\boldsymbol{U}^n$  defined on $\mathcal{F}(t^{n})$ to the full domain $\mathcal{O}$ by imposing the velocity on each node of $\mathcal{S}(t^{n+1})$ to be equal to ${\bf h}'^{n+1} + \theta'^{n+1}({\bf x}-{\bf h}^{n+1})^{\bot}$.\\
After this step, we know the Dirichlet condition for the velocity to impose at the interface
$\Gamma(t^{n+1}) = \partial \mathcal{S}(t^{n+1})$. So we determine $\boldsymbol{G}^{n+1}$ in~\eqref{M3} from ${\bf u}_{\Gamma}^{n+1} = {\bf h}'^{n+1} + \theta'^{n+1} ({\bf x}-{\bf h}^{n+1})^{\bot}$.
\item[5--] Finally, we compute $(\boldsymbol{U}^{n+1},\boldsymbol{P}^{n+1},\boldsymbol{\Lambda}^{n+1})$ such that
\begin{eqnarray*}
M_{{\bf u}{\bf u}} \frac{\boldsymbol{U}^{n+1} -\boldsymbol{U}^{n} }{dt} +A_{{\bf u}{\bf u}} \boldsymbol{U}^{n+1} + N(\boldsymbol{U}^{n+1}) \boldsymbol{U}^{n+1} +  A_{{\bf u}p}\boldsymbol{P}^{n+1}
+  A_{{\bf u}{\bf \lambda}}\boldsymbol{\Lambda}^{n+1} = \boldsymbol{F}^{n+1}, \\
A^T_{{\bf u}p}  \boldsymbol{U}^{n+1} +  A_{pp} \boldsymbol{P}^{n+1}  + A_{p{\bf \lambda}}\boldsymbol{\Lambda}^{n+1} = 0,  \\
A^T_{{\bf u}{\bf \lambda}}  \boldsymbol{U}^{n+1}  +  A^T_{p{\bf \lambda}} \boldsymbol{P}^{n+1}  +  A_{{\bf \lambda} {\bf \lambda}}\boldsymbol{\Lambda}^{n+1} = \boldsymbol{G}^{n+1}.
\end{eqnarray*}
At this stage, the solution of the resulting nonlinear algebraic system is achieved by a Newton method. The initialization of the Newton algorithm is done with the solution at the previous time step (this solution is defined in item 4--).
\item[6--] We complete the velocity $\boldsymbol{U}^{n+1}$ defined on $\mathcal{F}(t^{n+1})$ to the full domain $\mathcal{O}$ by imposing the velocity on each node of $\mathcal{S}(t^{n+1})$ to be equal to ${\bf h}'^{n+1} + \theta'^{n+1}({\bf x}-{\bf h}^{n+1})^{\bot}$.
\end{itemize}

{\it Remark 1:} In practice the mid-point method is used to update the geometry of the structure for the computation of ${\bf h}'^{n+1},{\bf h}^{n+1},{\theta'}^{n+1},\theta^{n+1}$.\\
{\it Remark 2:} Step 6-- plays an important role to update the geometry. Indeed, after extension, we have access to the values of the solution at each node of the full domain. Thus new nodes that appear after update have some values and no interpolation is required.

\section{Numerical tests: Free fall of a disk in a channel} \label{section4}
For validation of our method, we consider the numerical simulation of the motion of a disk falling inside an incompressible Newtonian viscous fluid. The parameters used in the computation, for a disk of radius $R=0.125$ cm in a channel of dimension $[0,2]\times [0,6]$, are given in Table \ref{tab}.
\begin{table}[h]
\begin{center}
\begin{tabular}{lllll}
Parameter& $\rho_f$ & $\rho_s$ & $\nu$ & $g$\\
\hline
Unit  & g/cm$^2$ & g/cm$^2$ & g/cm$^2$ s & cm/s$^2$\\  
value & 1 & 1.25& 0.1 & 981\\
\end{tabular}
\caption{Parameters used for the simulation.}
\label{tab}
\end{center}
\end{table}

This simulation is well documented in the literature and considered as a challenging benchmark.
We refer to the paper \cite{Glowinski} where fictitious domain method is used and \cite{Hachem} for simulations with mesh adaptation.\\
For the finite element discretization, we consider classical Lagrange family with $P_2 - P_1 - P_0$  for respectively ${\bf u}$, $p$, and $\boldsymbol{\lambda}$ , which is a choice that satisfies 
the inf-sup condition required for such kind of problems (see \cite{CourtFournieLozinski} for more details). Uniform triangular meshes are used and defined by imposing a uniform repartition of points on the boundary of the domain. Two meshes are used, $mesh_{50 \times 150}$ with $50$ points in $x$-direction and $150$ points in $y$-direction, and $mesh_{100 \times 300}$ with $100$ points in $x$-direction and $300$ points in $y$-direction.\\
Numerical tests are performed with and without stabilization to underline the advantage of the method. When stabilization is considered, we choose $\gamma =h\times \gamma_0$ where $\gamma_0=0.05$ (see  \cite{CourtFournieLozinski} for the justification of this choice). The parameter $\gamma$ has to obey to a compromise between the coerciveness of the system and the weight of the stabilization term. The time discretization step $dt$ is initialized to $0.0005$ and adapted at each time iteration to satisfy a CFL condition. More precisely, we evaluate the norm of the velocity at each point of the structure and we deduce the maximum value $v_m = \max_{{\bf x} \in \mathcal{S}(t^n)}(\|\bu({\bf x}) \|)$. Then we impose $dt = \min(0.9h/v_m,2h^2/{\nu})$. This condition is not restrictive and we observe that $dt \in [0.0005; 0.006]$ in all the tests.\\
In the literature, in order to study the fall of the disk, curves are given to show the evolution of the vertical velocity and the position of the center of the disk according to the time. We present the same analysis for different adaptation of our method. As expected, the disk reaches quickly a uniform fall velocity with slight moving on the right side of the vertical symmetry axis. This observation was already reported in the literature and is not specific to our method. One challenge is to propose robust methods that limit this breaking. In this section, we show that our method gives an answer to this question. Indeed, numerical simulations can be done with coarse meshes, even if it is not recommended with a fictitious domain approach (points on the interface can be far from the degrees of freedom introduced by the finite element method).\\

\begin{center}
\begin{figure}[h!]
\hspace*{3cm} \begin{minipage}{21cm}
\includegraphics[trim = 0cm 0cm 0cm 2cm, clip, scale=0.5]{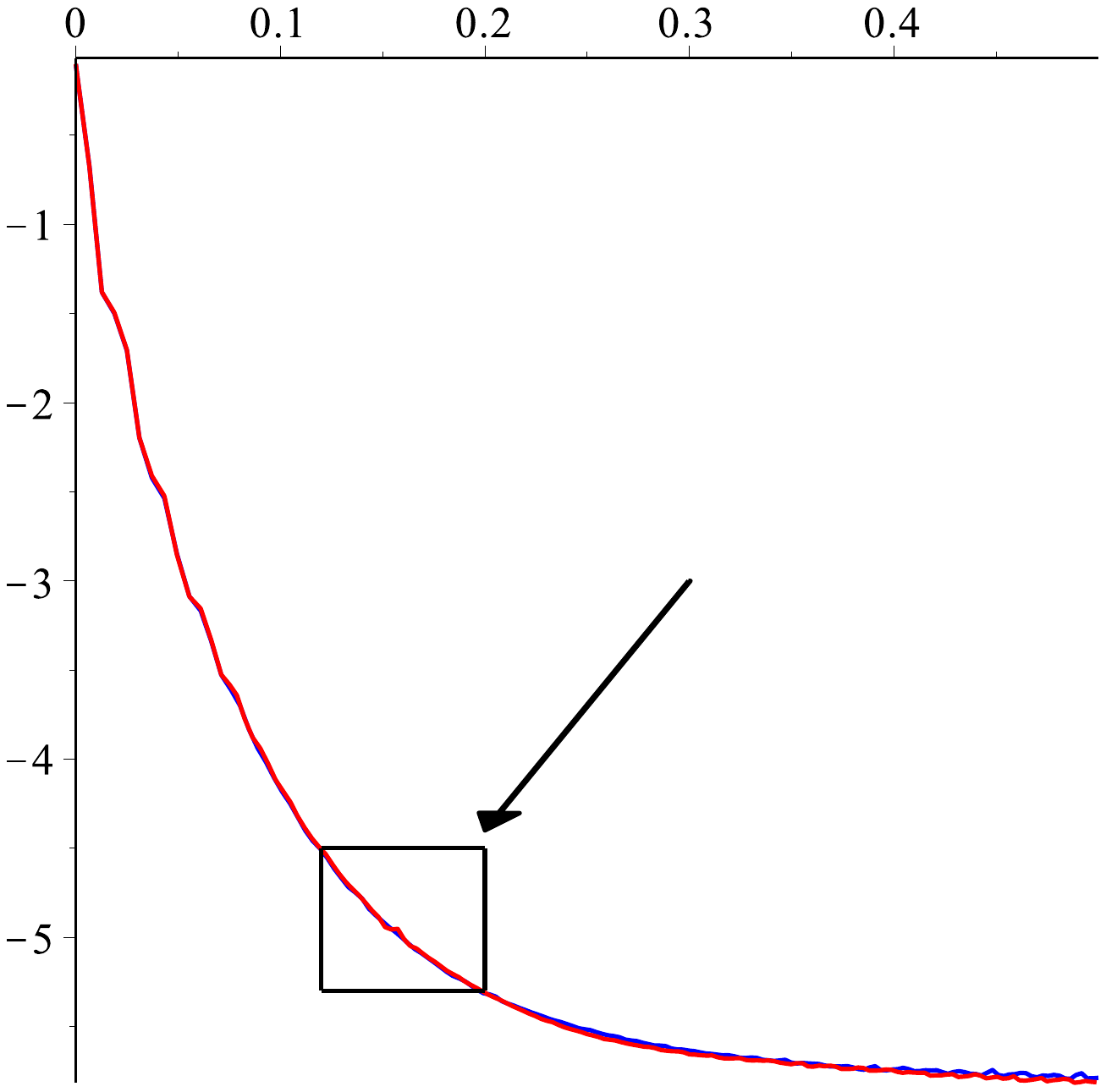} 

\vspace*{-11cm}
\hspace*{5.5cm}
\includegraphics[trim = 2cm 2cm 2cm 2cm, clip, scale=0.2]{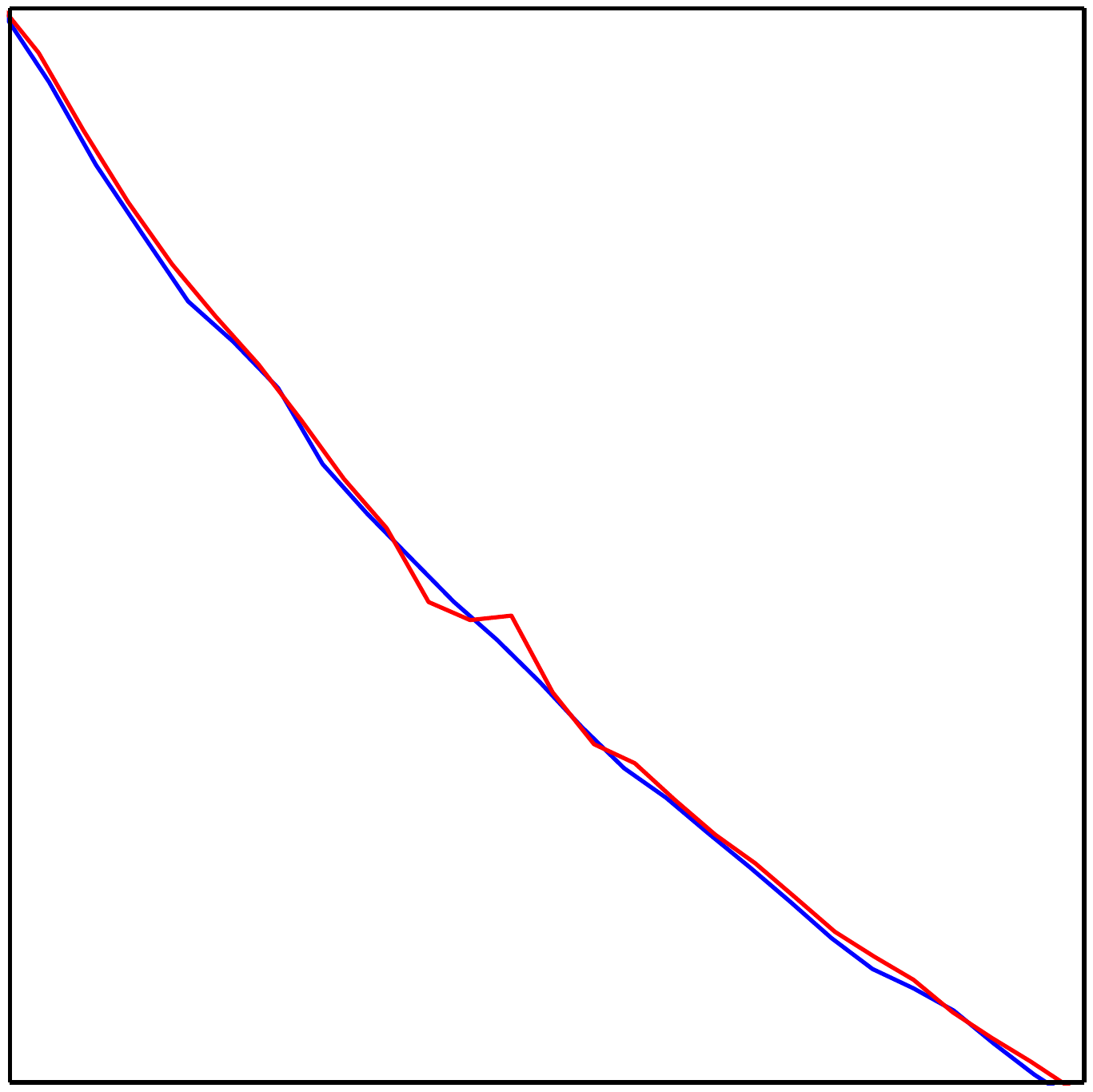}
\end{minipage}
\centering\caption{\label{FigGamma} Evolution in time $t\in[0;0.5]$ of the vertical velocity for $mesh_{50 \times 150}$ with stabilization in blue color and without stabilization in red color.}
\end{figure}
\end{center}

\begin{center}
\begin{figure}[h!]
\hspace*{-2cm} \begin{minipage}{20cm}
\hspace*{0cm} \includegraphics[trim = 0cm 0cm 5cm 0cm, clip, scale=0.5]{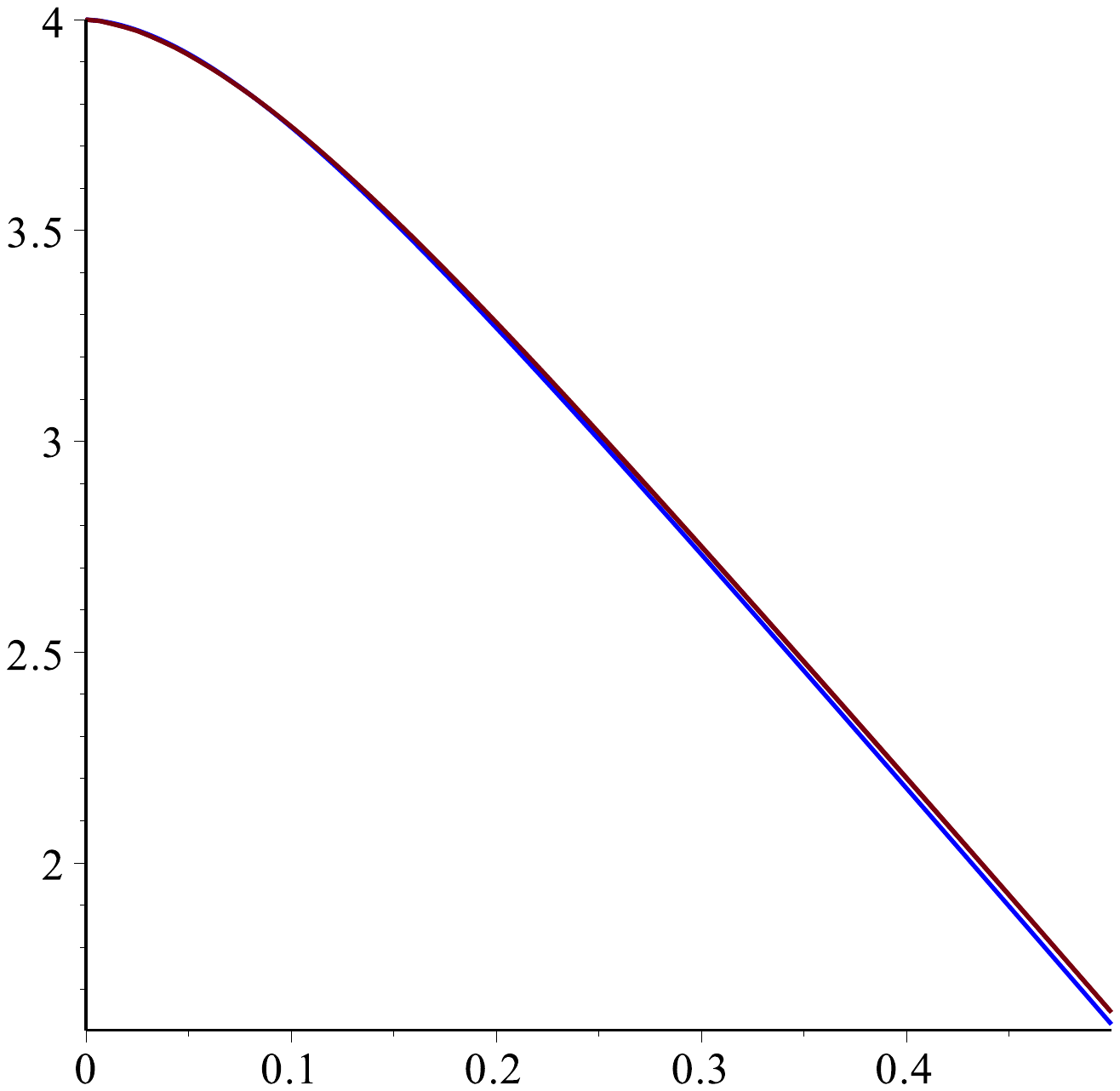} 
\vspace*{-14cm} 

\hspace*{8cm} \includegraphics[trim = 0cm 0cm 5cm 0cm, clip, scale=0.4]{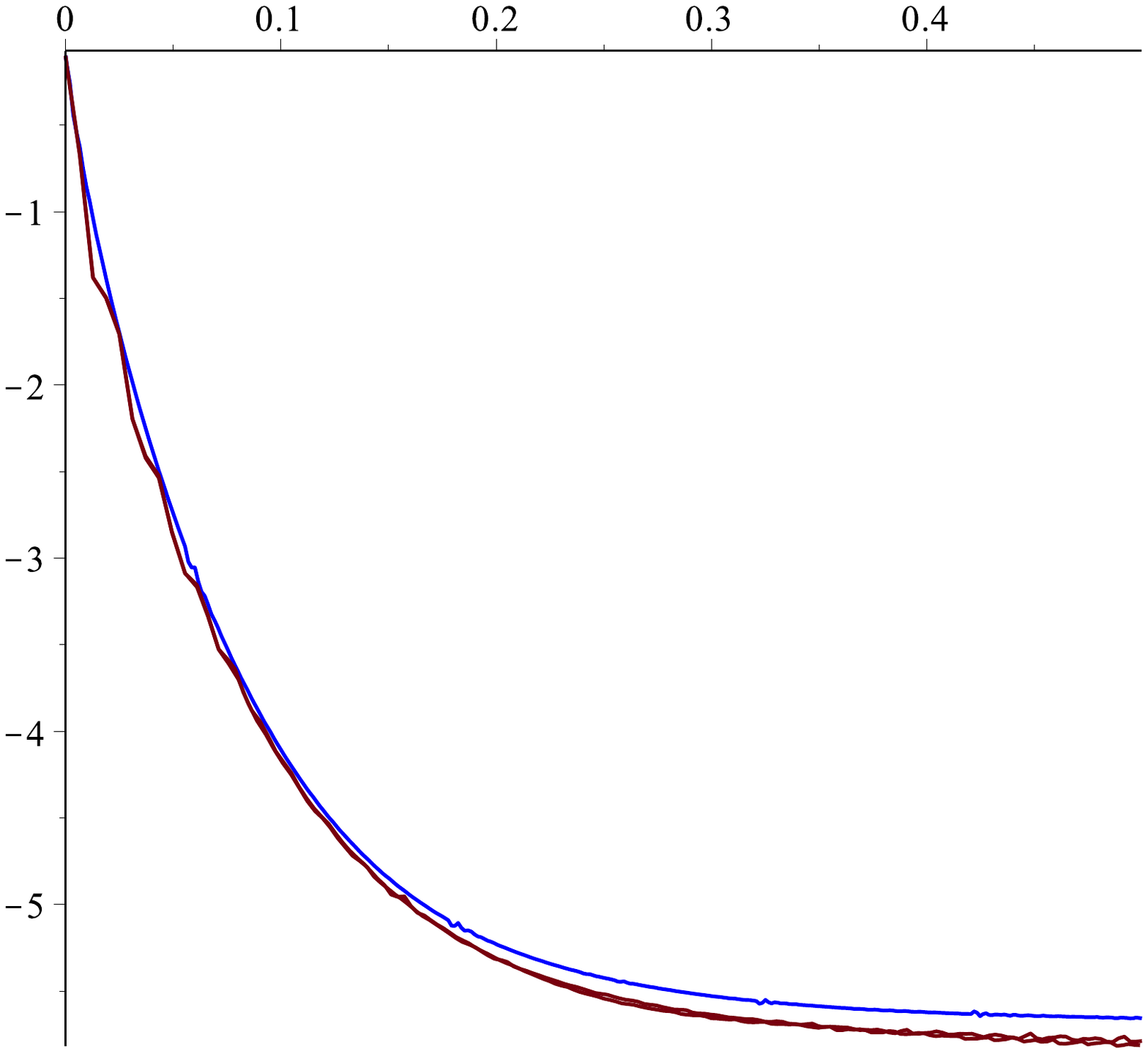} 
\end{minipage}
\vspace*{-3.5cm}

\hspace*{5cm} (a) \hspace*{7cm} (b)

\centering\caption{\label{FigCenterVelocity} Simulations for $t\in[0;0.5]$, for $mesh_{50 \times 150}$ in red color, and for $mesh_{100 \times 300}$ in blue color.
(a)~Vertical position of the disk,
(b)~Vertical velocity. }
\end{figure}
\end{center}

\begin{center}
\begin{figure}[h!]
\begin{center}
\hspace*{-2cm} \begin{minipage}{20cm}
 \includegraphics[trim = 0cm 7cm 0cm 0cm, clip, scale=0.5]{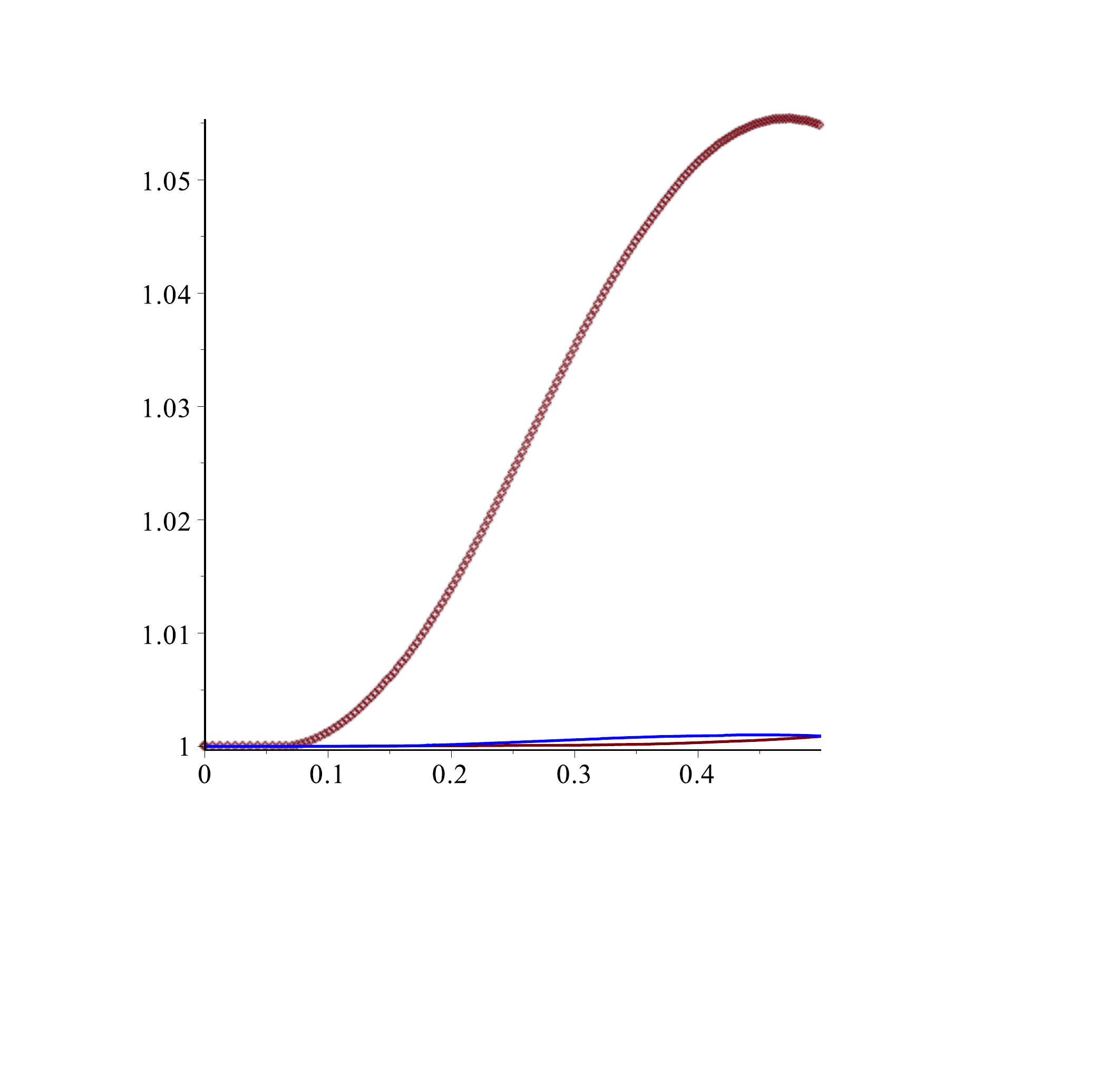} 
\vspace*{-10.5cm} 

\hspace*{8cm}
  \includegraphics[trim = 0cm 7cm 0cm 0cm, clip, scale=0.5]{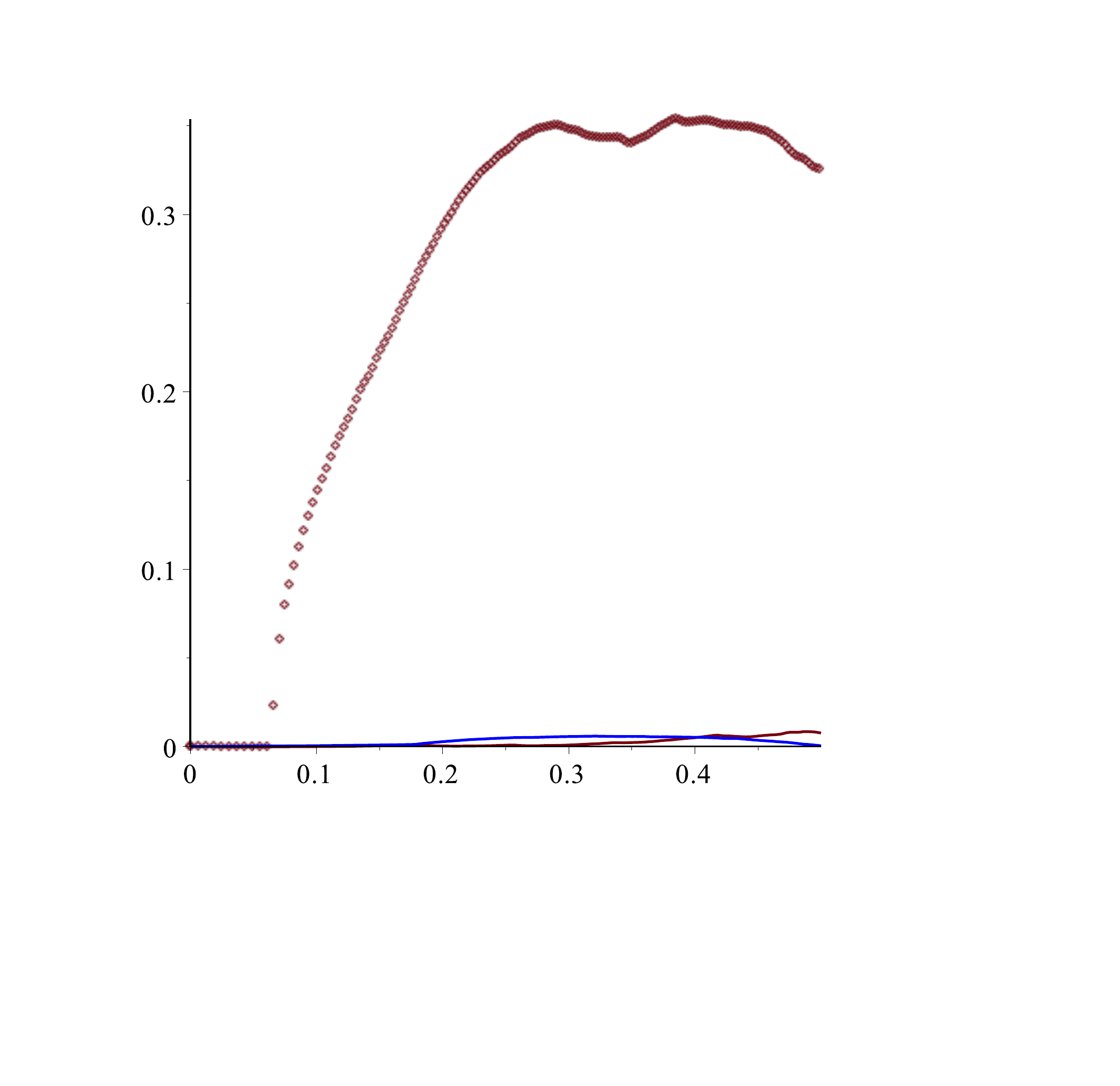} 
  \end{minipage}
\vspace*{-3cm}

\hspace*{0cm} (a) \hspace*{7cm} (b)

\caption{\label{FigPosx-angle} Simulations depending on time $t\in[0;0.5]$ with stabilization for $mesh_{50 \times 150}$  red lines and for $mesh_{100 \times 300}$ blue lines. 
Red curves (upper) with points have no stabilization with $mesh_{50 \times 150}$ points.
(a)~Horizontal position, 
(b)~Rotation angle.}
\end{center}
\end{figure}
\end{center}

{\bf Contribution of the stabilization technique}: Numerical tests are performed with the mesh $mesh_{50 \times 150}$, with and without performing the stabilization ($\gamma_0=0$). We compute the position of the disk according to the time and represent separately the vertical and horizontal positions.\\
In Figure~\ref{FigGamma} we represent the vertical velocity through the time. The results are similar, whether we perform stabilization or not. However, with the stabilization technique the method is more robust. If we zoom in (see Figure~\ref{FigGamma}), without stabilization (red curve) some perturbations appear. When we compare the positions of the disk through the time, we do not observe difference on the vertical position in the Figure~\ref{FigCenterVelocity}(a). However the difference is more important for the horizontal position and the rotation of the disk. The curves are plotted in red in Figure~\ref{FigPosx-angle}, with the mesh $mesh_{50\times 150}$. With the stabilization technique the results are clearly improved. Without stabilization, the symmetry is broken even if it seems that the disk comes back around the symmetry axis of the cavity at the end of the simulation. This behavior related to the computation of the angular velocity which is represented in Figure~\ref{FigPosx-angle}(b).

\FloatBarrier

{\bf Influence of the mesh size:} With stabilization, we compare the simulations obtained with $mesh_{50 \times 150}$ and $mesh_{100 \times 300}$. The results are given in Figure~\ref{FigCenterVelocity} (red curves for $mesh_{50 \times 150}$ and blue curves for $mesh_{100 \times 300}$). As expected, the smoothness of the solution is better when a sharper mesh is used. Besides, the result seems to be as good as results obtained in \cite{Glowinski, Hachem} for instance. When a coarse mesh is used, the velocity is over-estimated. This observation can be justified by the capability of the method for preserving the conservation of the mass. Indeed, with a coarse mesh, a numerical added mass appears in the system. This artificial mass is proportional to the stabilization terms (see $\mathcal{A}_{p{\blambda}} $  and $\mathcal{A}_{{\blambda}{\blambda}}$ in the discrete problem) which are themselves proportional to the mesh size, and thus it can be neglected when the mesh size decreases. 

The computation of the horizontal position of the disk is given in Figure~\ref{FigPosx-angle}(a). We observe that the disk tends to come back towards the symmetry axis of the cavity atthe end of the simulation with $mesh_{100 \times 300}$, unlike in the simulations with $mesh_{50 \times 150}$, where the symmetry breaking seems to growth. This phenomena can be observed in Figure~\ref{FigPosx-angle}(b) which represents the evolution of the rotation angle of the disk. With $mesh_{50 \times 150}$ this angle always growths, unlike for the sharper mesh $mesh_{100 \times 300}$. This behavior can be justified by perturbation associated with our numerical method, in particular for the treatment of the nonlinear term. Moreover, a perturbation in horizontal velocity component is difficult to compensate during the simulation and contributes to the amplifying the phenomena. However, with stabilization, relevant values of the rotation are obtained and show that their influence is reduced compared with the translation.\\

\FloatBarrier

\begin{center}
\begin{figure}[h!]
\begin{tabular} {ccccc}
\includegraphics[trim = 0cm 0cm 0cm 0cm, clip, scale=0.35]{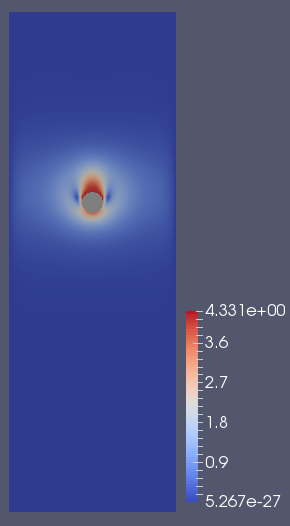} &
\hspace*{-10pt}
\includegraphics[trim = 0cm 0cm 0cm 0cm, clip, scale=0.35]{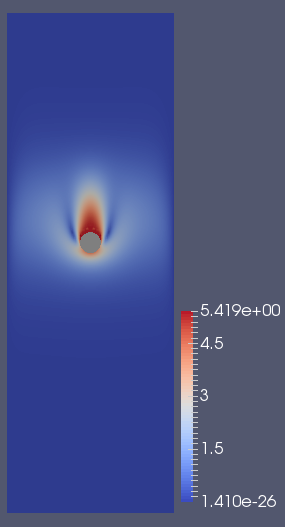} &
\hspace*{-10pt}
\includegraphics[trim = 0cm 0cm 0cm 0cm, clip, scale=0.35]{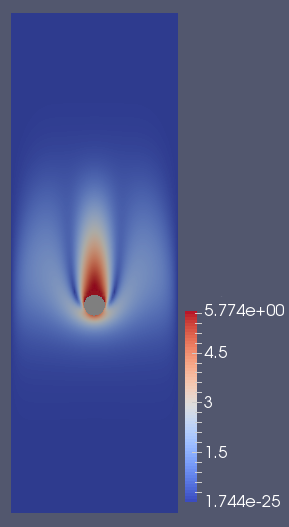} &
\hspace*{-10pt}
\includegraphics[trim = 0cm 0cm 0cm 0cm, clip, scale=0.35]{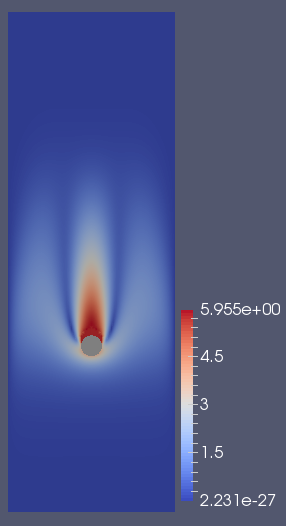} &
\hspace*{-10pt}
\includegraphics[trim = 0cm 0cm 0cm 0cm, clip, scale=0.35]{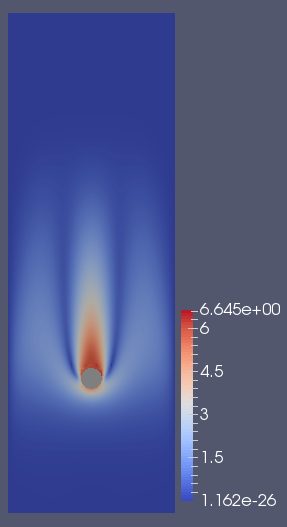}
\end{tabular}
\centering \caption{Imagery illustration of the intensity of the fluid's velocity during the fall of the ball.\label{figsuperfall}}
\end{figure}
\end{center}
\FloatBarrier

\section{Practical remarks on the numerical implementation.}
\begin{itemize}
\item All the numerical simulations were performed with the free generic library Getfem++ \cite{Getfem} (same source code for 2D and 3D) and implemented on High Performing Computers (parallel computations).
\item In order to compute properly the integrals over elements at the interface (during assembling procedure), external call to {\sc Qhull} Library \cite{qhull} is realized.
\item In the algorithm, steps 1-- and 2-- require computing of $M_{ \bf \lambda} \boldsymbol{\Lambda}^{n}$ and $M_{ \bf \lambda}  \left [ ({\bf x}-{\bf h}^{n})^{\bot}\cdot  \boldsymbol{\Lambda}^{n}\right ]$, corresponding to integrations over the level-set. Such integrations require particular attention, in the sake of preserving a good accuracy. Indeed, the integrations must use nodes on level-set and accurate values on that nodes are required (no interpolation).
\item The method is very efficient in time computation, since it requires an update of the assembling matrices only locally near the interface.
\item As mentioned in \cite{HaslR}, it is possible to define a reinforced stabilization technique in order to prevent difficulties that can occur when the intersection of the solid and the mesh over the whole domain introduce "very small" elements. The technique consists in selecting elements which are better to deduce the normal derivative on $\Gamma$. A similar approach is given in \cite{Pitkaranta}. We think that this kind of reinforced stabilization technique can prevent the perturbations that appear
during simulation (see zoom in Figure~\ref{FigGamma}).
\end{itemize}


\section{Conclusion}
In this paper, we have considered a new fictitious domain method based on the extended finite element with stabilized term applied to the Navier-Stokes equations coupled with a moving solid. This method is quite simple to implement since all the variables (multipliers and primal variables) are defined on a single mesh independent of the computational domain. The algorithm leads to a robust method (good computation of the normal Cauchy stress tensor) whatever is the intersection of the domain with the - not necessarily sharp - mesh. The simulation of a falling disk with respect to the time confirms that our approach is able to predict well the interaction between the fluid and the structure. The stabilization must be considered to obtain more physical results preserving symmetry. Applications in 3D are in progress, in particular for control flow by acting on the boundary of the solid.

\section*{Acknowledgements}
This work is partially supported by the foundation STAE in the context of the RTRA platform DYNAMORPH. It is based on the collaborative efforts with Alexei Lozinski and Yves Renard. \\


\small
\bibliographystyle{elsarticle-num}
\bibliography{court-fournie_Ref}

\begin{thebibliography}{10}
\expandafter\ifx\csname url\endcsname\relax
  \def\url#1{\texttt{#1}}\fi
\expandafter\ifx\csname urlprefix\endcsname\relax\def\urlprefix{URL }\fi
\expandafter\ifx\csname href\endcsname\relax
  \def\href#1#2{#2} \def\path#1{#1}\fi

\bibitem{Hou&Wang&Layton}
G.~Hou, J.~Wang, A.~Layton,
  \href{http://dx.doi.org/10.4208/cicp.291210.290411s}{Numerical methods for
  fluid-structure interaction---a review}, Commun. Comput. Phys. 12~(2) (2012)
  337--377.
\newblock \href {http://dx.doi.org/10.4208/cicp.291210.290411s}
  {\path{doi:10.4208/cicp.291210.290411s}}.
\newline\urlprefix\url{http://dx.doi.org/10.4208/cicp.291210.290411s}

\bibitem{LT}
G.~Legendre, T.~Takahashi,
  \href{http://dx.doi.org/10.1051/m2an:2008020}{Convergence of a
  {L}agrange-{G}alerkin method for a fluid-rigid body system in {ALE}
  formulation}, M2AN Math. Model. Numer. Anal. 42~(4) (2008) 609--644.
\newblock \href {http://dx.doi.org/10.1051/m2an:2008020}
  {\path{doi:10.1051/m2an:2008020}}.
\newline\urlprefix\url{http://dx.doi.org/10.1051/m2an:2008020}

\bibitem{SMSTT0}
J.~San~Mart{\'{\i}}n, J.-F. Scheid, T.~Takahashi, M.~Tucsnak,
  \href{http://dx.doi.org/10.1137/S0036142903438161}{Convergence of the
  {L}agrange-{G}alerkin method for the equations modelling the motion of a
  fluid-rigid system}, SIAM J. Numer. Anal. 43~(4) (2005) 1536--1571.
\newblock \href {http://dx.doi.org/10.1137/S0036142903438161}
  {\path{doi:10.1137/S0036142903438161}}.
\newline\urlprefix\url{http://dx.doi.org/10.1137/S0036142903438161}

\bibitem{SST}
J.~San~Mart{\'{\i}}n, L.~Smaranda, T.~Takahashi,
  \href{http://dx.doi.org/10.1016/j.cam.2008.12.021}{Convergence of a finite
  element/{ALE} method for the {S}tokes equations in a domain depending on
  time}, J. Comput. Appl. Math. 230~(2) (2009) 521--545.
\newblock \href {http://dx.doi.org/10.1016/j.cam.2008.12.021}
  {\path{doi:10.1016/j.cam.2008.12.021}}.
\newline\urlprefix\url{http://dx.doi.org/10.1016/j.cam.2008.12.021}

\bibitem{Peskinacta}
C.~S. Peskin, \href{http://dx.doi.org/10.1017/S0962492902000077}{The immersed
  boundary method}, Acta Numer. 11 (2002) 479--517.
\newblock \href {http://dx.doi.org/10.1017/S0962492902000077}
  {\path{doi:10.1017/S0962492902000077}}.
\newline\urlprefix\url{http://dx.doi.org/10.1017/S0962492902000077}

\bibitem{Mittal&Iaccarino}
R.~Mittal, G.~Iaccarino,
  \href{http://dx.doi.org/10.1146/annurev.fluid.37.061903.175743}{Immersed
  boundary methods}, in: Annual review of fluid mechanics. {V}ol. 37, Vol.~37
  of Annu. Rev. Fluid Mech., Annual Reviews, Palo Alto, CA, 2005, pp. 239--261.
\newblock \href {http://dx.doi.org/10.1146/annurev.fluid.37.061903.175743}
  {\path{doi:10.1146/annurev.fluid.37.061903.175743}}.
\newline\urlprefix\url{http://dx.doi.org/10.1146/annurev.fluid.37.061903.175743}

\bibitem{Glowinski}
R.~Glowinski, T.-W. Pan, T.~I. Hesla, D.~D. Joseph, A distributed lagrange
  multiplier / fictitious domain method for particular flows, Int. J. of
  Multiphase Flow 25 (1999) 755--794.

\bibitem{MoesD}
N.~Mo{\"e}s, J.~Dolbow, T.~Belytschko, A finite element method for crack growth
  without remeshing, Internat. J. Numer. Methods Engrg. 46~(1) (1999) 131--150.
\newblock \href
  {http://dx.doi.org/10.1002/(SICI)1097-0207(19990910)46:1<131::AID-NME726>3.0.CO;2-J}
  {\path{doi:10.1002/(SICI)1097-0207(19990910)46:1<131::AID-NME726>3.0.CO;2-J}}.

\bibitem{reviewXfem}
T.-P. Fries, T.~Belytschko, \href{http://dx.doi.org/10.1002/nme.2914}{The
  extended/generalized finite element method: an overview of the method and its
  applications}, Internat. J. Numer. Methods Engrg. 84~(3) (2010) 253--304.
\newblock \href {http://dx.doi.org/10.1002/nme.2914}
  {\path{doi:10.1002/nme.2914}}.
\newline\urlprefix\url{http://dx.doi.org/10.1002/nme.2914}

\bibitem{MoesB}
N.~Mo{\"e}s, E.~B{\'e}chet, M.~Tourbier,
  \href{http://dx.doi.org/10.1002/nme.1675}{Imposing {D}irichlet boundary
  conditions in the extended finite element method}, Internat. J. Numer.
  Methods Engrg. 67~(12) (2006) 1641--1669.
\newblock \href {http://dx.doi.org/10.1002/nme.1675}
  {\path{doi:10.1002/nme.1675}}.
\newline\urlprefix\url{http://dx.doi.org/10.1002/nme.1675}

\bibitem{SukumarC}
N.~Sukumar, D.~L. Chopp, N.~Mo{\"e}s, T.~Belytschko,
  \href{http://dx.doi.org/10.1016/S0045-7825(01)00215-8}{Modeling holes and
  inclusions by level sets in the extended finite-element method}, Comput.
  Methods Appl. Mech. Engrg. 190~(46-47) (2001) 6183--6200.
\newblock \href {http://dx.doi.org/10.1016/S0045-7825(01)00215-8}
  {\path{doi:10.1016/S0045-7825(01)00215-8}}.
\newline\urlprefix\url{http://dx.doi.org/10.1016/S0045-7825(01)00215-8}

\bibitem{Gerstenberger2008}
A.~Gerstenberger, W.~A. Wall,
  \href{http://dx.doi.org/10.1016/j.cma.2007.07.002}{An extended finite element
  method/{L}agrange multiplier based approach for fluid-structure interaction},
  Comput. Methods Appl. Mech. Engrg. 197~(19-20) (2008) 1699--1714.
\newblock \href {http://dx.doi.org/10.1016/j.cma.2007.07.002}
  {\path{doi:10.1016/j.cma.2007.07.002}}.
\newline\urlprefix\url{http://dx.doi.org/10.1016/j.cma.2007.07.002}

\bibitem{Choi2010}
Y.~J. Choi, M.~A. Hulsen, H.~E.~H. Meijer,
  \href{http://dx.doi.org/doi:10.1016/j.jnnfm.2010.02.021}{An extended finite
  element method for the simulation of particulate viscoelastic flows}, J.
  Non-Newtonian Fluid Mech. 165~(11-12) (2010) 607--624.
\newblock \href {http://dx.doi.org/doi:10.1016/j.jnnfm.2010.02.021}
  {\path{doi:doi:10.1016/j.jnnfm.2010.02.021}}.
\newline\urlprefix\url{http://dx.doi.org/doi:10.1016/j.jnnfm.2010.02.021}

\bibitem{HaslR}
J.~Haslinger, Y.~Renard, \href{http://dx.doi.org/10.1137/070704435}{A new
  fictitious domain approach inspired by the extended finite element method},
  SIAM J. Numer. Anal. 47~(2) (2009) 1474--1499.
\newblock \href {http://dx.doi.org/10.1137/070704435}
  {\path{doi:10.1137/070704435}}.
\newline\urlprefix\url{http://dx.doi.org/10.1137/070704435}

\bibitem{MoesG}
N.~Mo{\"e}s, A.~Gravouil, T.~Belytschko,
  \href{http://dx.doi.org/10.1002/nme.429}{Non-planar 3d crack growth by the
  extended finite element and level sets, part i: Mechanical model}, Internat.
  J. Numer. Methods Engrg. 53~(11) (2002) 2549--2568.
\newblock \href {http://dx.doi.org/10.1002/nme.429}
  {\path{doi:10.1002/nme.429}}.
\newline\urlprefix\url{http://dx.doi.org/10.1002/nme.429}

\bibitem{Stazi}
F.~L. Stazi, E.~Budyn, J.~Chessa, T.~Belytschko, An extended finite element
  method with high-order elements for curved cracks, Comput. Mech. 31~(1-2)
  (2003) 38--48.

\bibitem{SukumarM}
N.~Sukumar, N.~Mo\"{e}s, B.~Moran, T.~Belytschko, Extended finite element
  method for three-dimensional crack modelling, Int. J. Numer. Meth. Engng
  48~(11) (2000) 1549--1570.

\bibitem{Stolarska}
M.~Stolarska, D.~L. Chopp, N.~Mo\"{e}s, T.~Belytschko,
  \href{http://dx.doi.org/10.1002/nme.201}{Modelling crack growth by level sets
  in the extended finite element method}, Int. J. Numer. Meth. Engng 51~(8)
  (2001) 943--960.
\newblock \href {http://dx.doi.org/10.1002/nme.201}
  {\path{doi:10.1002/nme.201}}.
\newline\urlprefix\url{http://dx.doi.org/10.1002/nme.201}

\bibitem{Bechet2009}
{\'E}.~B{\'e}chet, N.~Mo{\"e}s, B.~Wohlmuth,
  \href{http://dx.doi.org/10.1002/nme.2515}{A stable {L}agrange multiplier
  space for stiff interface conditions within the extended finite element
  method}, Internat. J. Numer. Methods Engrg. 78~(8) (2009) 931--954.
\newblock \href {http://dx.doi.org/10.1002/nme.2515}
  {\path{doi:10.1002/nme.2515}}.
\newline\urlprefix\url{http://dx.doi.org/10.1002/nme.2515}

\bibitem{HildR}
P.~Hild, Y.~Renard, \href{http://dx.doi.org/10.1007/s00211-009-0273-z}{A
  stabilized {L}agrange multiplier method for the finite element approximation
  of contact problems in elastostatics}, Numer. Math. 115~(1) (2010) 101--129.
\newblock \href {http://dx.doi.org/10.1007/s00211-009-0273-z}
  {\path{doi:10.1007/s00211-009-0273-z}}.
\newline\urlprefix\url{http://dx.doi.org/10.1007/s00211-009-0273-z}

\bibitem{CourtFournieLozinski}
S.~Court, M.~Fourni{\'e}, A.~Lozinski,
  \href{http://dx.doi.org/10.1002/fld.3839}{A fictitious domain approach for
  the {S}tokes problem based on the extended finite element method}, Internat.
  J. Numer. Methods Fluids 74~(2) (2014) 73--99.
\newblock \href {http://dx.doi.org/10.1002/fld.3839}
  {\path{doi:10.1002/fld.3839}}.
\newline\urlprefix\url{http://dx.doi.org/10.1002/fld.3839}

\bibitem{Nitsche}
J.~Nitsche, \"{U}ber ein {V}ariationsprinzip zur {L}\"osung von
  {D}irichlet-{P}roblemen bei {V}erwendung von {T}eilr\"aumen, die keinen
  {R}andbedingungen unterworfen sind, Abh. Math. Sem. Univ. Hamburg 36 (1971)
  9--15, collection of articles dedicated to Lothar Collatz on his sixtieth
  birthday.

\bibitem{Burman1}
R.~Becker, E.~Burman, P.~Hansbo, \href{http://dx.doi.org/10.1002/nme.3093}{A
  hierarchical {NXFEM} for fictitious domain simulations}, Internat. J. Numer.
  Methods Engrg. 86~(4-5) (2011) 549--559.
\newblock \href {http://dx.doi.org/10.1002/nme.3093}
  {\path{doi:10.1002/nme.3093}}.
\newline\urlprefix\url{http://dx.doi.org/10.1002/nme.3093}

\bibitem{Burman3}
E.~Burman, P.~Hansbo,
  \href{http://dx.doi.org/10.1016/j.apnum.2011.01.008}{Fictitious domain finite
  element methods using cut elements: {II}. {A} stabilized {N}itsche method},
  Appl. Numer. Math. 62~(4) (2012) 328--341.
\newblock \href {http://dx.doi.org/10.1016/j.apnum.2011.01.008}
  {\path{doi:10.1016/j.apnum.2011.01.008}}.
\newline\urlprefix\url{http://dx.doi.org/10.1016/j.apnum.2011.01.008}

\bibitem{Massing}
A.~Massing, M.~G. Larson, A.~Logg, M.~E. Rognes,
  \href{http://dx.doi.org/10.1007/s10915-014-9838-9}{A stabilized {N}itsche
  fictitious domain method for the {S}tokes problem}, J. Sci. Comput. 61~(3)
  (2014) 604--628.
\newblock \href {http://dx.doi.org/10.1007/s10915-014-9838-9}
  {\path{doi:10.1007/s10915-014-9838-9}}.
\newline\urlprefix\url{http://dx.doi.org/10.1007/s10915-014-9838-9}

\bibitem{JPR}
C.~Airiau, J.~M. Buchot, R.~K. Dubey, M.~Fourni\'e, J.~P. Raymond, Best
  actuator location for stabilization of the navier-stokes equations,
  Submitted.

\bibitem{Evans}
L.~C. Evans, Partial differential equations, 2nd Edition, Vol.~19 of Graduate
  Studies in Mathematics, American Mathematical Society, Providence, RI, 2010.

\bibitem{Gunzburger}
M.~D. Gunzburger, S.~L. Hou, \href{http://dx.doi.org/10.1137/0729024}{Treating
  inhomogeneous essential boundary conditions in finite element methods and the
  calculation of boundary stresses}, SIAM J. Numer. Anal. 29~(2) (1992)
  390--424.
\newblock \href {http://dx.doi.org/10.1137/0729024}
  {\path{doi:10.1137/0729024}}.
\newline\urlprefix\url{http://dx.doi.org/10.1137/0729024}

\bibitem{Barbosa1}
H.~J.~C. Barbosa, T.~J.~R. Hughes,
  \href{http://dx.doi.org/10.1016/0045-7825(91)90125-P}{The finite element
  method with {L}agrange multipliers on the boundary: circumventing the
  {B}abu\v ska-{B}rezzi condition}, Comput. Methods Appl. Mech. Engrg. 85~(1)
  (1991) 109--128.
\newblock \href {http://dx.doi.org/10.1016/0045-7825(91)90125-P}
  {\path{doi:10.1016/0045-7825(91)90125-P}}.
\newline\urlprefix\url{http://dx.doi.org/10.1016/0045-7825(91)90125-P}

\bibitem{Ern}
A.~Ern, J.-L. Guermond,
  \href{http://dx.doi.org/10.1007/978-1-4757-4355-5}{Theory and practice of
  finite elements}, Vol. 159 of Applied Mathematical Sciences, Springer-Verlag,
  New York, 2004.
\newblock \href {http://dx.doi.org/10.1007/978-1-4757-4355-5}
  {\path{doi:10.1007/978-1-4757-4355-5}}.
\newline\urlprefix\url{http://dx.doi.org/10.1007/978-1-4757-4355-5}

\bibitem{Hachem}
E.~Hachem, S.~Feghali, R.~Codina, T.~Coupez, Anisotropic adaptive meshing and
  monolithic variational multiscale method for fluid–structure interaction.

\bibitem{Getfem}
Y.~Renard, J.~Pommier, \href{http://home.gna.org/getfem/}{An open source
  generic C++ library for finite element methods}.
\newline\urlprefix\url{http://home.gna.org/getfem/}

\bibitem{qhull}
C.~B. Barber, D.~P. Dobkin, H.~Huhdanpaa,
  \href{http://dx.doi.org/10.1145/235815.235821}{The quickhull algorithm for
  convex hulls}, ACM Trans. Math. Software 22~(4) (1996) 469--483.
\newblock \href {http://dx.doi.org/10.1145/235815.235821}
  {\path{doi:10.1145/235815.235821}}.
\newline\urlprefix\url{http://dx.doi.org/10.1145/235815.235821}

\bibitem{Pitkaranta}
J.~Pitk{\"a}ranta, \href{http://dx.doi.org/10.2307/2006378}{Local stability
  conditions for the {B}abu\v ska method of {L}agrange multipliers}, Math.
  Comp. 35~(152) (1980) 1113--1129.
\newblock \href {http://dx.doi.org/10.2307/2006378}
  {\path{doi:10.2307/2006378}}.
\newline\urlprefix\url{http://dx.doi.org/10.2307/2006378}

\end{thebibliography}

\end{document}